\documentclass[11pt]{amsart}
\usepackage{hyperref}
\usepackage[margin=1in]{geometry}
 \pdfoutput=1 
\usepackage{amssymb}
\usepackage{epstopdf}
\usepackage{float}
\usepackage{amsmath}

\usepackage{tikz}
\usetikzlibrary{decorations.pathreplacing,calligraphy}
\usetikzlibrary {turtle} 
\usetikzlibrary{math}
\usetikzlibrary{positioning,calc}
\usetikzlibrary{arrows}
\usetikzlibrary{decorations.markings}
\usepackage[ruled, lined, linesnumbered, commentsnumbered, longend]{algorithm2e}

\newcommand\comment[1]{}

\def\mylist{{0}{90}{180}{270}}
\pgfmathdeclarerandomlist{myang}{\mylist}

\begin{document}

\title{L-systems for the Boundaries of fractal dragon
  space filling folding curves}

  \author{{
H. A. Verrill
\date{\today}
H.A.Verrill@warwick.ac.uk
    } 
    }

    \maketitle

\newtheorem{definition}{Definition}
\newtheorem{notation}{Notation}
\newtheorem{corollary}{Corollary}
\newtheorem{lemma}{Lemma}
\newtheorem{example}{Example}
\newtheorem{theorem}{Theorem}
\newtheorem{proposition}{Proposition}
\newtheorem{assumption}{Assumption}

\newcommand\diagtiles[8]{
        \begin{scope}[shift = {(#1+#3/2,#2+#3/2)}]
\begin{scope}[scale={#7}]
          \begin{scope}[rotate={#6}]
          \draw[gray, fill=#5](-#3/2,-#3/2) rectangle (#3/2,#3/2);
          \begin{scope}[rotate = {#4}]
            \draw[black, line width = 3] (1/2,#3/2)--(-1/2,-#3/2);
            \draw[black, line width = 3,->] (1/2,#3/2)--(0,0);            
    \end{scope}
    \end{scope}    
        \end{scope}
      \end{scope}
}

\newcommand\R{{\texttt{R}}}
\renewcommand\L{{\texttt{L}}}
\renewcommand\r{{\texttt{r}}}
\renewcommand\l{{\texttt{l}}}
\renewcommand\S{{\texttt{S}}}
\newcommand\s{{\texttt{s}}}
\newcommand\A{{\texttt{A}}}
\newcommand\B{{\texttt{B}}}
\newcommand\+{{\texttt{+}}}
\renewcommand{\-}{\text{\texttt{-}}}
\newcommand{\myminus}{\text{\texttt{-}}}
\renewcommand\v{{\texttt{v}}}

\newcommand\planetiles[8]{
        \begin{scope}[shift = {(#1+#3/2,#2+#3/2)}]
\begin{scope}[scale={#7}]
          \begin{scope}[rotate={#6}]
          \draw[gray, fill=#5](-#3/2,-#3/2) rectangle (#3/2,#3/2);
          \begin{scope}[rotate = {#4}]
    \end{scope}
    \end{scope}    
        \end{scope}
      \end{scope}
}

             \newcommand\tileCCC[8]{
        \begin{scope}[shift = {(#1+#3/2,#2+#3/2)}]
\begin{scope}[scale={#7}]
          \begin{scope}[rotate={#6}]
          \draw[gray, fill=#5](-#3/2,-#3/2) rectangle (#3/2,#3/2);
    \begin{scope}[rotate = {#4}]
      \draw[<-,line width = #8,blue] (-#3/2,0) arc (-90:0:#3/2);
\draw[->,line width = #8,blue] (0,-#3/2) arc (180:90:#3/2);    
    \end{scope}
    \end{scope}    
        \end{scope}
      \end{scope}
            }

             \newcommand\tileCCCx[8]{
        \begin{scope}[shift = {(#1+#3/2,#2+#3/2)}]
\begin{scope}[scale={#7}]
          \begin{scope}[rotate={#6}]
          \draw[gray, fill=#5](-#3/2,-#3/2) rectangle (#3/2,#3/2);
    \begin{scope}[rotate = {#4}]
      \draw[<-,line width = #8] (-#3/2,0) arc (-90:0:#3/2);
\draw[->,line width = #8] (0,-#3/2) arc (180:90:#3/2);    
    \end{scope}
    \end{scope}    
        \end{scope}
      \end{scope}
            }

            \newcommand\tileDDD[8]{
        \begin{scope}[shift = {(#1+#3/2,#2+#3/2)}]
\begin{scope}[scale={#7}]
          \begin{scope}[rotate={#6}]
          \draw[gray, fill=#5](-#3/2,-#3/2) rectangle (#3/2,#3/2);
    \begin{scope}[rotate = {#4}]
      \draw[->,line width = #8,blue] (-#3/2,0) arc (-90:0:#3/2);
\draw[<-,line width = #8,blue] (0,-#3/2) arc (180:90:#3/2);    
    \end{scope}
    \end{scope}    
        \end{scope}
      \end{scope}
            }

             \newcommand\tileCCCb[8]{
        \begin{scope}[shift = {(#1+#3/2,#2+#3/2)}]
\begin{scope}[scale={#7}]
          \begin{scope}[rotate={#6}]
          \draw[gray, fill=#5](-#3/2,-#3/2) rectangle (#3/2,#3/2);
    \begin{scope}[rotate = {#4}]
      \draw[<-,line width = #8] (-#3/2,0) arc (-90:0:#3/2);
      \draw[->,line width = #8] (0,-#3/2) arc (180:90:#3/2);
          \draw[orange, line width = 3] (1/2,#3/2)--(-1/2,-#3/2);
    \end{scope}
    \end{scope}    
        \end{scope}
      \end{scope}
            }

            \newcommand\tileDDDb[8]{
        \begin{scope}[shift = {(#1+#3/2,#2+#3/2)}]
\begin{scope}[scale={#7}]
          \begin{scope}[rotate={#6}]
          \draw[gray, fill=#5](-#3/2,-#3/2) rectangle (#3/2,#3/2);
    \begin{scope}[rotate = {#4}]
      \draw[->,line width = #8,blue] (-#3/2,0) arc (-90:0:#3/2);
      \draw[<-,line width = #8,blue
      ] (0,-#3/2) arc (180:90:#3/2);
    \draw[line width = 3] (1/2,#3/2)--(-1/2,-#3/2);     
    \end{scope}
    \end{scope}    
        \end{scope}
      \end{scope}
            }

                         \newcommand\tileCCCbb[8]{
        \begin{scope}[shift = {(#1+#3/2,#2+#3/2)}]
\begin{scope}[scale={#7}]
          \begin{scope}[rotate={#6}]
          \draw[gray, fill=#5](-#3/2,-#3/2) rectangle (#3/2,#3/2);
    \begin{scope}[rotate = {#4}]
          \draw[orange, line width = 3] (1/2,#3/2)--(-1/2,-#3/2);
    \end{scope}
    \end{scope}    
        \end{scope}
      \end{scope}
            }

            \newcommand\tileDDDbb[8]{
        \begin{scope}[shift = {(#1+#3/2,#2+#3/2)}]
\begin{scope}[scale={#7}]
          \begin{scope}[rotate={#6}]
          \draw[gray, fill=#5](-#3/2,-#3/2) rectangle (#3/2,#3/2);
    \begin{scope}[rotate = {#4}]
    \draw[orange, line width = 3] (1/2,#3/2)--(-1/2,-#3/2);     
    \end{scope}
    \end{scope}    
        \end{scope}
      \end{scope}
            }

            \newcommand\tileDDDbbb[8]{
        \begin{scope}[shift = {(#1+#3/2,#2+#3/2)}]
\begin{scope}[scale={#7}]
          \begin{scope}[rotate={#6}]
          \draw[gray, fill=#5](-#3/2,-#3/2) rectangle (#3/2,#3/2);
    \begin{scope}[rotate = {#4}]
      \draw[->,line width = #8] (-#3/2,0) arc (-90:0:#3/2);
      \draw[<-,line width = #8] (0,-#3/2) arc (180:90:#3/2);
      \draw[orange, line width = 3] (1/2,#3/2)--(-1/2,-#3/2);
          \draw[orange, line width = 3,->] (1/2,#3/2)--(0,0);     
    \end{scope}
    \end{scope}    
        \end{scope}
      \end{scope}
  }               

            \newcommand\tileAB[5]{
  \begin{scope}[shift = {(#1+#3/2,#2+#3/2)}]
    \begin{scope}[rotate = {#4}]
    \draw[gray,line width = 2] (-#3/2,0) arc (-90:0:#3/2);
\draw[gray,line width = 2] (0,-#3/2) arc (180:90:#3/2);    
    \end{scope}
      \end{scope}
            }

            \newcommand\tileABb[5]{
  \begin{scope}[shift = {(#1+#3/2,#2+#3/2)}]
    \begin{scope}[rotate = {#4}]
    \draw[gray,line width = 2] (-#3/2,0) arc (-90:0:#3/2);
    \draw[gray,line width = 2] (0,-#3/2) arc (180:90:#3/2);
    \draw[orange, line width = 3] (1/2,#3/2)--(-1/2,-#3/2);
    \end{scope}
      \end{scope}
  }

\newcommand\twistBB[1]{
         \clip (0,0) rectangle (1,1);
          \begin{scope}[shift = {(0.5,0.5)}]
          \begin{scope}[scale = {0.86},rotate={-10}]
          \begin{scope}[shift = {(-0.5,-0.5)}]
            \tileAB{1}{0}{1}{0}{gray!50!white};
            \tileAB{0}{1}{1}{90}{gray!50!white};
            \tileAB{0}{-1}{1}{90}{gray!50!white};
            \tileAB{-1}{0}{1}{0}{gray!50!white};            
         \end{scope}
          \end{scope}
                  \end{scope}
}

\newcommand\twistBBb[1]{
         \clip (0,0) rectangle (1,1);
          \begin{scope}[shift = {(0.5,0.5)}]
          \begin{scope}[scale = {0.86},rotate={-10}]
          \begin{scope}[shift = {(-0.5,-0.5)}]
            \tileABb{1}{0}{1}{0}{gray!50!white};
            \tileABb{0}{1}{1}{90}{gray!50!white};
            \tileABb{0}{-1}{1}{90}{gray!50!white};
            \tileABb{-1}{0}{1}{0}{gray!50!white};            
         \end{scope}
          \end{scope}
                  \end{scope}
          }

\newcommand\twistAA[1]{
         \clip (0,0) rectangle (1,1);
          \begin{scope}[shift = {(0.5,0.5)}]
          \begin{scope}[scale = {0.86},rotate={10}]
          \begin{scope}[shift = {(-0.5,-0.5)}]
            \tileAB{0}{0}{1}{#1}{gray!50!white};
            \tileAB{1}{0}{1}{90}{gray!50!white};
            \tileAB{0}{1}{1}{0}{gray!50!white};
            \tileAB{0}{-1}{1}{0}{gray!50!white};
            \tileAB{-1}{0}{1}{90}{gray!50!white};            
         \end{scope}
          \end{scope}
          \end{scope}          
}

\newcommand\twistAAb[1]{
         \clip (0,0) rectangle (1,1);
          \begin{scope}[shift = {(0.5,0.5)}]
          \begin{scope}[scale = {0.86},rotate={10}]
          \begin{scope}[shift = {(-0.5,-0.5)}]
            \tileABb{0}{0}{1}{#1}{gray!50!white};
            \tileABb{1}{0}{1}{90}{gray!50!white};
            \tileABb{0}{1}{1}{0}{gray!50!white};
            \tileABb{0}{-1}{1}{0}{gray!50!white};
            \tileABb{-1}{0}{1}{90}{gray!50!white};            
         \end{scope}
          \end{scope}
          \end{scope}          
}

\newcommand\twistBC[1]{
         \clip (0,0) rectangle (1,1);
          \begin{scope}[shift = {(0.5,0.5)}]
          \begin{scope}[scale = {0.71},rotate={-45}]
          \begin{scope}[shift = {(-0.5,-0.5)}]
            \tileAB{1}{0}{1}{0}{gray!50!white};
            \tileAB{0}{1}{1}{90}{gray!50!white};
            \tileAB{0}{-1}{1}{90}{gray!50!white};
            \tileAB{-1}{0}{1}{0}{gray!50!white};            
         \end{scope}
          \end{scope}
                  \end{scope}
}

\newcommand\twistBCb[1]{
         \clip (0,0) rectangle (1,1);
          \begin{scope}[shift = {(0.5,0.5)}]
          \begin{scope}[scale = {0.71},rotate={-45}]
          \begin{scope}[shift = {(-0.5,-0.5)}]
            \tileABb{1}{0}{1}{0}{gray!50!white};
            \tileABb{0}{1}{1}{90}{gray!50!white};
            \tileABb{0}{-1}{1}{90}{gray!50!white};
            \tileABb{-1}{0}{1}{0}{gray!50!white};            
         \end{scope}
          \end{scope}
                  \end{scope}
          }

\newcommand\twistAC[1]{
         \clip (0,0) rectangle (1,1);
          \begin{scope}[shift = {(0.5,0.5)}]
          \begin{scope}[scale = {0.71},rotate={45}]
          \begin{scope}[shift = {(-0.5,-0.5)}]
            \tileAB{0}{0}{1}{#1}{gray!50!white};
            \tileAB{1}{0}{1}{90}{gray!50!white};
            \tileAB{0}{1}{1}{0}{gray!50!white};
            \tileAB{0}{-1}{1}{0}{gray!50!white};
            \tileAB{-1}{0}{1}{90}{gray!50!white};            
         \end{scope}
          \end{scope}
          \end{scope}          
}

\newcommand\mypath{
   (0,0)--++(0,1)--++(1,0)
    --++(0,1)
    --++(-1,0)
    --++(0,1)
    --++(1,0)
    --++(0,-1)
    --++(1,0)
    --++(0,-1)
    --++(-1,0)
    --++(0,-1)
    --++(1,0)
    --++(0,1)
    --++(1,0)
    --++(0,1)
        --++(1,0)
        --++(0,-1);
        }

\newcommand\mypathA{
   (0,0)--++(1,0)
    --++(0,1)
    --++(1,0)
    --++(0,1)
    --++(1,0)
    --++(0,-1)
    --++(1,0)
    --++(0,-1)
    --++(-1,0);
        }

   \section*{Abstract}
    We describe an algorithm to find an
    L-system for the boundary of plane-filling square grid based folding curves, such as the
    fractal dragon curves.

\section{Introduction: plane-filling folding curves}
Our starting point is Heighway's fractal dragon \cite{Edgar}, and its extensions, for example, as described in
\cite{arndt}, \cite{dekking}, \cite{Tab}, which should be consulted for further
details of the plane-filling curve construction and notation.
An iterative sequence of paths is constructed in terms of an L-system.
The limit is a plane-filling curve.
Examples are shown in Figure~\ref{fig:3examples}.
\begin{figure}
  \begin{tikzpicture}
    \begin{scope}[shift={(-8,1)}]
    \draw[thin, gray](2.5,-1.5) rectangle (6.5,2.5);
    \clip(2.5,-1.5) rectangle (6.5,2.5);
  \foreach\i in {0,...,4}{
    \foreach\j in {0,...,4}{
      \draw[line width=1mm,->,red] (\i+\j+0.1,\i-\j)--(\i+\j+0.9,\i-\j);
      \draw[line width=1mm,->,red] (\i+\j-0.1,\i-\j)--(\i+\j-0.9,\i-\j);
      \draw[line width=1mm,->,blue] (\i+\j+1,\i-\j+0.1)--(\i+\j+1,\i-\j+0.9);
      \draw[line width=1mm,->,blue] (\i+\j+1,\i-\j-0.1)--(\i+\j+1,\i-\j-0.9);      
      \node at (\i+\j+0.5,\i-\j+0.15){\tiny\texttt{A}};
      \node at (\i+\j-0.5,\i-\j+0.15){\tiny\texttt{A}};
      \node at (\i+\j+0.15,\i-\j+0.5){\tiny\texttt{B}};
      \node at (\i+\j+0.15,\i-\j-0.5){\tiny\texttt{B}};
  }}
  \end{scope}

            \tileCCC{0}{0}{1}{0}{gray!50!white}{0}{1}{2}
            \tileCCC{2}{0}{1}{0}{gray!50!white}{0}{1}{2}            
            \tileDDD{0}{2}{1}{90}{gray!50!white}{0}{1}{2}
            \tileDDD{1}{1}{1}{90}{gray!50!white}{0}{1}{2}
            \tileDDD{2}{2}{1}{90}{gray!50!white}{0}{1}{2}

            \tileDDD{2}{1}{1}{90}{white}{90}{1}{2}
            \tileDDD{1}{2}{1}{90}{white}{90}{1}{2}
            \tileDDD{0}{1}{1}{90}{white}{90}{1}{2}
            \tileCCC{1}{0}{1}{0}{white}{90}{1}{2}

\draw[line width=1mm,red](1,0)++(135:0.5) arc (135:45:0.5);
\draw[line width=1mm,red](1,0)++(135:0.5) arc (135:45:0.5);
\draw[line width=1mm,red](1,1)++(135:0.5) arc (135:45:0.5);
\draw[line width=1mm,red](3,0)++(90:0.5) arc (90:135:0.5);
\draw[line width=1mm,red](3,1)++(90:0.5) arc (90:135:0.5);
\draw[line width=1mm,red](0,2)++(90:0.5) arc (90:45:0.5);
\draw[line width=1mm,red](0,1)++(-90:0.5) arc (-90:-45:0.5);
\draw[line width=1mm,red](3,3)++(-90:0.5) arc (-90:-135:0.5);
\draw[line width=1mm,red](0,2)++(-90:0.5) arc (-90:-45:0.5);
\draw[line width=1mm,red](2,2)++(135:0.5) arc (135:45:0.5);
\draw[line width=1mm,red](2,2)++(-135:0.5) arc (-135:-45:0.5);
\draw[line width=1mm,red](2,1)++(-135:0.5) arc (-135:-45:0.5);
\draw[line width=1mm,red](1,3)++(-135:0.5) arc (-135:-45:0.5);

            \node at (0.5,0.5){$\+$};
            \node at (0.5,1.5){$\myminus$};
            \node at (0.5,2.5){$\myminus$};
            \node at (1.5,0.5){$\+$};
            \node at (1.5,1.5){$\myminus$};
            \node at (1.5,2.5){$\myminus$};
            \node at (2.5,0.5){$\+$};
            \node at (2.5,1.5){$\myminus$};
            \node at (2.5,2.5){$\myminus$};

            \begin{scope}[shift={(4,0)}]

\input{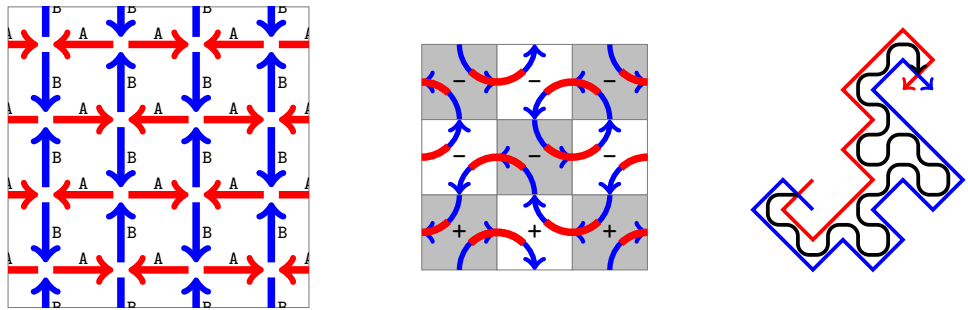}
              
\end{scope}

  \end{tikzpicture}
  \caption{
    Left: Square grid of allowed path segments. Middle: Choice of turns, $+$ or $-$ drawn on dual grid, with curved segments.
    Right: Heighway dragon and its boundary after 5 iterations
    of the L-systems.}
  \label{fig:grid1}
  \end{figure}
Paths of the plane-filling curve are composed of a sequence of alternating \texttt{A}
and \texttt{B} edges, which lie on
a square grid as in Figure~\ref{fig:grid1} (a).  At each vertex, a choice left, indicated by
\texttt{-}, or right, indicated by \texttt{+}, must be made.
This is summarised in Table~\ref{tab:ABL}
\begin{table}
  \begin{tabular}{|lll|lll|}
    \hline
         \texttt{A}&: & move forward one unit, starting from an even square &  \texttt{+}&: & turn right\\
         \texttt{B}&: & move forward one unit, starting from an odd square &        \texttt{-}&: &turn left\\
    \hline
  \end{tabular}
     \caption{Geometric interpretation of the alphabet of the L-system for the folding curves.
     }
     \label{tab:ABL}
      \end{table}
Paths are not allowed to cross or overlap, so at each intersection, both incoming
paths must turn the same way.  In Figure~\ref{fig:grid1} (b), the paths are redrawn with a choice of \texttt{+}
or \texttt{-} 
labeled.  These segments are drawn with rounded corners,
which allows us to see the direction of turning more clearly.
These paths are now drawn on a dual lattice, with vertices replaced by squares.
We colour squares in a checkerboard pattern, of gray and white squares,
which we call even and odd respectively. 

An L-system consists of a triple, $(\Omega,\mathcal A,P)$, where
$\Omega$ is an alphabet, $\mathcal A$ is a word in $\Omega$ called an axiom, and
$P$ is a function which maps elements of $\Omega$ to words in $\Omega$.
We follow the notation of \cite{A}, \cite{arndt} and \cite{dekking}.
Our plane-filling curves correspond to an
alphabet $\Omega:=\{\texttt{A},\texttt{B},\texttt{+},\texttt{-}\}$, and words have the form
\begin{equation}
  X_1 s_1 X_2 s_2 X_3s_3 \dots s_{n-1}X_n, \text{ where }
   X_i\in\{\texttt{A},\texttt{B}\}, X_{i}\not=X_{i+1}, s_i\in\{\+,\-\}
\label{eqn:1}
\end{equation}
The L-system rule has the form of a folding curve rule \cite{dekking}:
$$
\begin{array}{lll}
\texttt{A}\mapsto \sigma(\texttt{A})&:=& X_0 s_1 X_1 s_2 X_2 \dots s_n X_n\\
\texttt{B}\mapsto \sigma(\texttt{B})&:=&\tilde X_n \tilde s_n \tilde X_{n-1} \tilde s_{n-1} \tilde X_{n-2} \dots \tilde s_1 \tilde X_0,\\
\texttt{+}\mapsto \sigma(\texttt{+})&:=&\texttt{+}\\
\texttt{-}\mapsto \sigma(\texttt{-})&:=&\texttt{-}\\
\end{array}
$$

where $\tilde{\texttt{A}} = \texttt{B}$,
$\tilde{\texttt{B}}=\texttt{A}$,
$\tilde{\texttt{+}}=\texttt{-}, \tilde{\texttt{-}}=\texttt{+}$.

That is, $\sigma(\texttt{B})$ is obtained by reversing $\sigma(\texttt{A})$ and switching $\texttt{A}\leftrightarrow \texttt{B}$
and \texttt{+}$\leftrightarrow$ \texttt{-}.
Further requirements are necessary for this system to be the L-system of a plane-filling curve, see \cite{dekking} and \cite{arndt}.
For the Heighway fractal dragon, the rule  $\sigma$ is as in Table~\ref{tab:fill}, case 1.
Repeated application of this rule results in the Heighway fractal dragon curve.  Application of $\sigma$ five
times to a segment {\A}
is shown in Figure~\ref{fig:grid1}(c), black curve.
There are many variants \cite{Tab}.
This paper is concerned with the boundaries of these folding curves.

\section{The L-system for the boundary of plane-filling curves}

We describe an algorithm, which given a square grid folding curve L-system, as
described in \cite{dekking} and \cite{arndt},
given by an L-system rule $\sigma$, will produce an L-system $\tau$ for the boundary of the curve.
More precisely, the curves
$\tau^n(\texttt{R})$ and $\tau^n(\texttt{L})$ are approximations to the left and right boundaries of
$\sigma^n(\texttt{A})$, and the boundary of $\lim_{n\rightarrow \infty}\sigma^n(\texttt{A})$ is
the union of  $\lim_{n\rightarrow \infty}\sigma^n(\R)$ and  $\lim_{n\rightarrow \infty}\sigma^n(\texttt{L})$.
Examples are given by the red and blue curves in Figures~\ref{fig:grid1}(c) and Figure~\ref{fig:3examples}.

The L-system $\tau$ has alphabet consisting of 6 symbols,
\texttt{L}, \texttt{R}, \texttt{l}, \texttt{r}, \texttt{S}, and \texttt{s}.
These are interpreted as in Table~\ref{tab:B}.

\begin{figure}
  \begin{tikzpicture}[scale={0.75}]

    \diagtiles{1}{0}{1}{90}{white!50!white}{0}{1}{2}
            \diagtiles{1}{1}{1}{-90}{gray!60!white}{90}{1}{2}
            \planetiles{0}{1}{1}{90}{white!50!white}{0}{1}{2}
            \planetiles{0}{0}{1}{90}{gray!50!white}{0}{1}{2}
\node at (0.75,1){\L};           
\node at (1,-0.5){\L};           

\begin{scope}[shift={(2.5,0))}]
    \planetiles{1}{0}{1}{90}{white!50!white}{0}{1}{2}
            \diagtiles{1}{1}{1}{90}{gray!60!white}{90}{1}{2}
            \diagtiles{0}{1}{1}{90}{white!50!white}{0}{1}{2}
            \planetiles{0}{0}{1}{90}{gray!50!white}{0}{1}{2}
\node at (1,0.75){\l};           
\node at (1,-0.5){\l};           

           \end{scope} 

\begin{scope}[shift={(5,0)}]
    \diagtiles{1}{0}{1}{-90}{white!50!white}{0}{1}{2}
            \diagtiles{1}{1}{1}{90}{gray!60!white}{90}{1}{2}
            \planetiles{0}{1}{1}{-90}{white!50!white}{0}{1}{2}
            \planetiles{0}{0}{1}{-90}{gray!50!white}{0}{1}{2}
\node at (0.75,1){\r};           
\node at (1,-0.5){\r};           

\begin{scope}[shift={(2.5,0))}]
    \planetiles{1}{0}{1}{-90}{white!50!white}{0}{1}{2}
            \diagtiles{1}{1}{1}{-90}{gray!60!white}{90}{1}{2}
            \diagtiles{0}{1}{1}{-90}{white!50!white}{0}{1}{2}
            \planetiles{0}{0}{1}{-90}{gray!50!white}{0}{1}{2}
\node at (1,0.75){\R};           
\node at (1,-0.5){\R};           

\end{scope}
           \end{scope} 

\begin{scope}[shift={(10,0)}]
    \planetiles{1}{0}{1}{-90}{white!50!white}{0}{1}{2}
            \diagtiles{1}{1}{1}{90}{gray!60!white}{90}{1}{2}
            \planetiles{0}{1}{1}{-90}{white!50!white}{0}{1}{2}
            \diagtiles{0}{0}{1}{180}{gray!50!white}{0}{1}{2}
\node at (0.75,1){\S};           
\node at (1,-0.5){\S};           

\begin{scope}[shift={(2.5,0))}]
    \diagtiles{1}{0}{1}{-90}{white!50!white}{0}{1}{2}
            \planetiles{1}{1}{1}{-90}{gray!60!white}{90}{1}{2}
            \diagtiles{0}{1}{1}{-90}{white!50!white}{0}{1}{2}
            \planetiles{0}{0}{1}{-90}{gray!50!white}{0}{1}{2}
\node at (1,0.75){\texttt{s}};           
\node at (1,-0.5){\s};           

\end{scope}
           \end{scope}

    \end{tikzpicture}
  \caption{Geometric interpretation of letters of L-system. Gray squares are even, white odd.   
    Starting  and ending    points  centre of squares.  
  }
  \label{fig:kindsofconnections}
\end{figure}
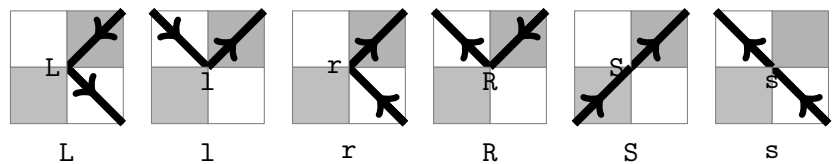

\begin{figure}
          \begin{tikzpicture}[scale={0.8}]

            \diagtiles{0}{-1}{1}{180}{gray!50!white}{0}{1}{2}
            \planetiles{1}{-1}{1}{90}{white!50!white}{0}{1}{2}
            \planetiles{2}{-1}{1}{90}{gray!50!white}{0}{1}{2}
            \diagtiles{3}{-1}{1}{90}{white!50!white}{0}{1}{2}
            \diagtiles{4}{-1}{1}{180}{gray!50!white}{0}{1}{2}

            \diagtiles{0}{0}{1}{-90}{white!50!white}{0}{1}{2}
            \diagtiles{0}{1}{1}{90}{gray!60!white}{90}{1}{2}
            \planetiles{0}{2}{1}{90}{white!50!white}{0}{1}{2}

            \planetiles{1}{0}{1}{90}{gray!60!white}{90}{1}{2}
            \planetiles{1}{1}{1}{-90}{white!50!white}{0}{1}{2}
            \diagtiles{1}{2}{1}{90}{gray!60!white}{90}{1}{2}
            
            \diagtiles{2}{0}{1}{90}{white!50!white}{0}{1}{2}            
            \diagtiles{2}{1}{1}{0}{gray!60!white}{0}{1}{2}
            \diagtiles{2}{2}{1}{90}{white!50!white}{0}{1}{2}
            \planetiles{3}{0}{1}{180}{gray!50!white}{0}{1}{2}            
            \planetiles{3}{1}{1}{180}{white!60!white}{90}{1}{2}
            \planetiles{3}{2}{1}{180}{gray!50!white}{0}{1}{2}

                        \diagtiles{4}{0}{1}{-90}{white!50!white}{0}{1}{2}            
            \diagtiles{4}{1}{1}{90}{gray!60!white}{90}{1}{2}
            \diagtiles{4}{2}{1}{-90}{white!50!white}{0}{1}{2}

            \node at (1,-0.25){\L};
            \node at (-0.25,1){\r};
            \node at (1,2.25){\S};
            \node at (2,3.25){\R};
            \node at (2.75,2){\r};
            \node at (1.75,1){\L};
            \node at (3,-0.25){\s};
            \node at (4,-1.25){\l};
            \node at (5.25,0){\L};                                              
            \node at (4,0.75){\r};
            \node at (5.25,2){\L};

\begin{scope}[shift={(6,0)},scale={2}]
\draw (0,0) grid (2,1);
\draw[fill=gray!60!white](0,0) rectangle (1,1);
\draw[fill=white!50!white](2,0) rectangle (1,1);
\filldraw (0.5,0.5) circle (0.1);
\filldraw (1.5,0.5) circle (0.1);

\draw[red,line width=1mm,->](0.5,0.5) --(1,0.5);
\draw[red,line width=1mm](1,.5) --(1.5,0.5);
\draw[line width=1mm,->](0.5,0.5)--(1,1);
\draw[line width=1mm](1,1)--(1.5,0.5);
\draw[line width=1mm,->](0.5,0.5)--(1,0);
\draw[line width=1mm](1,0)--(1.5,0.5);
\node at(1.2,0.6){\A};
\node at(1,1.2){\R};
\node at(1,-0.2){\L};
\end{scope}

  \begin{scope}[shift={(11,0)}]
       \tileDDDb{0}{0}{1}{90}{white!50!white}{0}{1}{2}
            \tileDDDb{2}{0}{1}{90}{white!50!white}{0}{1}{2}            
                        \tileDDDb{0}{2}{1}{90}{white!50!white}{0}{1}{2}
                        
            \tileDDDb{1}{1}{1}{90}{white!50!white}{0}{1}{2}
            \tileCCC{2}{2}{1}{0}{white!50!white}{0}{1}{2} 

            \tileCCC{2}{1}{1}{0}{gray!60!white}{90}{1}{2}
            \tileDDDb{1}{2}{1}{90}{gray!60!white}{90}{1}{2}
            \tileCCC{0}{1}{1}{0}{gray!60!white}{90}{1}{2}
            \tileDDDb{1}{0}{1}{90}{gray!60!white}{90}{1}{2}

          \end{scope}

          \end{tikzpicture}
          \caption{Left: Path corresponding to the word
            \texttt{LrSRrLslLrL}. Middle: relationship between {\A} and {\L} and \R,
            which are all paths from the mid point of an even square to the
            mid point of an odd square.
            Right: boundary segments (black), which go in diagonal directions
            on the same grid as $\A, \B$ segments (both blue),
            which go horizontally and vertically
            respectively. (Though $\A, \B$ are shown curved, so this is an approximation.)
          }
          \label{fig:newword}
\end{figure}
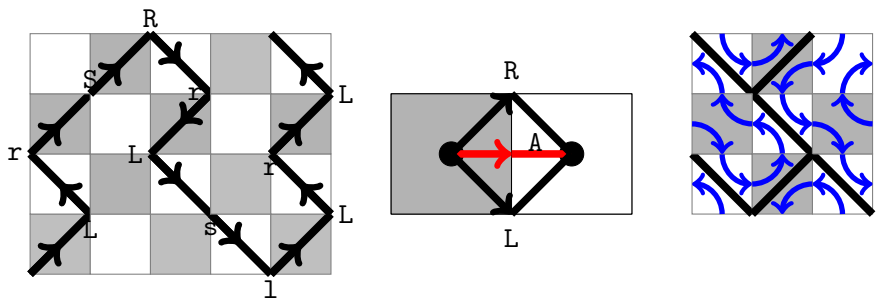

\begin{figure}
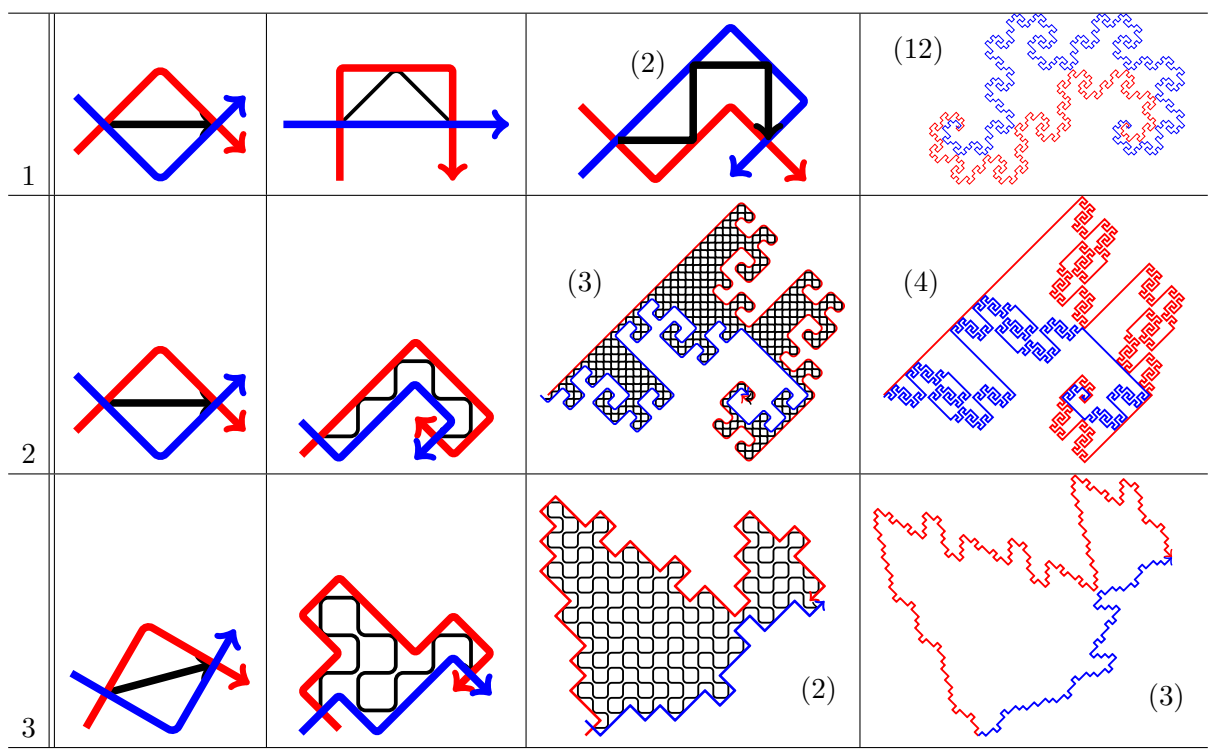

$$
\begin{array}{r||c|c|c|c}
    \hline
  1&
  \begin{tikzpicture}[scale={0.75}]
    \draw[line width = 1mm,->](0,0)--(2,0);
                    \draw[line width = 1mm,->,red,rounded corners=1mm](-0.5,-0.5)--++(1.5,1.5)--++(1.5,-1.5);
      \draw[line width = 1mm,->,blue,rounded corners=1mm](-0.5,0.5)--++(1.5,-1.5)--++(1.5,1.5);      
    \end{tikzpicture}
&
          \begin{tikzpicture}[scale={0.75}]
        \draw[black,line width = 0.5mm,rounded corners=1mm,->](0,0)--(1,1)--(2,0);
        \draw [->,rounded corners=0.5mm, line width = 1mm, red] (0,-1)--(0,1)--(2,1)--(2,-1);
          \draw [->,rounded corners=0.5mm, line width = 1mm, blue] (-1,0)--(3,0);            
          \end{tikzpicture}
          &
 \begin{tikzpicture}
  \draw [->,rounded corners=0.5mm, red, line width = 1mm,xshift={-4mm},yshift={4.6mm},  turtle/distance=1.412,turtle={home,right=135,
   forward,   left, forward, right, forward
}];

  \draw [->,rounded corners=0.5mm, blue, line width = 1mm, xshift={-4mm},yshift={-4.8mm},turtle/distance=1.412,turtle={home,right=45,
      forward, forward, right, forward, right, forward, 
}];

  \draw [->,xshift=1mm,rounded corners=0.25mm, line width = 1mm, turtle/distance=1,turtle={home,right=90,
      forward, left, forward, right, forward, right, forward
    }];  

    \node at (0.5,1){(2)};
      
\end{tikzpicture}          
 &
  \input{example2FF.tex}          
          \\
            \hline
  2&
  \begin{tikzpicture}[scale={0.75}]
    \draw[line width = 1mm,->](0,0)--(2,0);
          \draw[line width = 1mm,->,red,rounded corners=1mm](-0.5,-0.5)--++(1.5,1.5)--++(1.5,-1.5);
      \draw[line width = 1mm,->,blue,rounded corners=1mm](-0.5,0.5)--++(1.5,-1.5)--++(1.5,1.5);      
    \end{tikzpicture}
    &
          \begin{tikzpicture}
      \begin{scope}[scale={0.5},shift={(0,4)}]
\draw[black,line width = 0.5mm,rounded corners=1mm,->]\mypathA
    \end{scope}
    \input{example3BA.tex}
          \end{tikzpicture}
          &
          \input{exampleC2.tex}
          &
\input{exampleB2.tex}          
    \\
  \hline
  3 &
    \begin{tikzpicture}[scale={0.75},rotate={15}]
      \draw[line width = 1mm,->](0,0)--(2,0);
                \draw[line width = 1mm,->,red,rounded corners=1mm](-0.5,-0.5)--++(1.5,1.5)--++(1.5,-1.5);
      \draw[line width = 1mm,->,blue,rounded corners=1mm](-0.5,0.5)--++(1.5,-1.5)--++(1.5,1.5);      
    \end{tikzpicture}
&
      \begin{tikzpicture}
      \begin{scope}[scale={0.5},shift={(0,4)}]
\draw[black,line width = 0.5mm,rounded corners=1mm,->]\mypath
    \end{scope}
    \draw [->,rounded corners=0.5mm, line width = 1mm, red, yshift={1.75cm}, xshift={0.25cm},turtle/distance=7.15mm,turtle={home,left=45,
  forward,
right, forward, left, forward, right, forward, right, forward, forward, left, forward, right, forward, right,forward
}];

\draw [->,rounded corners=0.5mm, line width = 1mm,blue, yshift={1.7cm}, xshift={-0.26cm},turtle/distance=7.15mm,turtle={home,left=-45,forward,
right, forward, left, forward, forward, right,forward
}];
  \end{tikzpicture}
&
    \begin{tikzpicture}
     \input{example3.tex}
           \input{example3DD.tex}
    \end{tikzpicture}
    &
               \begin{tikzpicture}
           \input{example3E.tex}
               \end{tikzpicture}
               \\
                 \hline
\end{array}
$$
\caption{Three examples of plane-filling curves and boundaries, with L-systems as in
  Tables~\ref{tab:fill} and \ref{tab:bou}.  Column 1 shows segment {\A} (black), a right turn \R (red), on the left side of
  \A, and a left turn, {\L} (blue), on the right side of {\A}.
  Column 2 shows \A, \L, {\R}, after one application of the relevant L-system.
  The number of iterations for columns 3 and 4 is given in brackets.}
\label{fig:3examples}
\end{figure}

\begin{table}
     \begin{tabular}{lll}
         \texttt{L}&: & turn left, and move forward one unit, starting from an even square\\
         \texttt{R}&:  &turn right, and move forward one unit, starting from an even square\\ 
         \texttt{l}&: & turn left, and move forward one unit, starting from an odd square\\
       \texttt{r}&: &turn right, and move forward one unit, starting from an odd square\\ 
       \texttt{S}&: &do not turn;  move forward one unit, starting from an even square\\
       \texttt{s}&: &do not turn;  move forward one unit, starting from an odd square\\
     \end{tabular}
     \caption{The geometric interpretation of the L-system alphabet for the boundary of plane-filling folding curves.
     }
     \label{tab:B}
\end{table}
Figure~\ref{fig:kindsofconnections} shows the geometric interpretation of each case of
${\texttt{L, l, R, r, S, s}}$, extended to
cover two tiles.
These symbols may also be interpreted as ``move half a unit, follow indicated direction, move half a unit''.
This 
       results in the same path, up to slight difference at the ends.  Here ``unit'' means the whole of the diagonal of one
       square.   Change of direction on the diagonal, in the middle of a square is not permitted, except
       possibly at the beginning and end of the original path, for example, to change from $\tau(\texttt{R})$
       to $\tau(\texttt{L})$,
       which otherwise live on different sublattices to each other, that is, the turning vertices of $\tau(\R)$ never
       conincide with those of $\tau(\texttt{L})$.
Figure~\ref{fig:newword} (left) shows an example of how the word \texttt{LrSRrLslLrL} is interpreted as a path.
Note that we can divide the vertices of this grid into two sub-lattices, one of which has \texttt{L,R,S} lables,
the other has \texttt{l,r,s} labels.

\section{The Algorithm}

To compute the L-system rule $\tau$, start with the rule given by the sequence $\sigma(\A)$.
Note that we do not address the question of finding $\sigma$; we assume $\sigma$ gives an L-system for a
plane-filling folding
curve.  Examples of $\sigma$ are given the next section, and further examples
can be found for example, in \cite{dekking}, \cite{arndt}, \cite{AH}.
The algorithm is as follows.  Note that in addition to
\texttt{A,B,+,-,R,r,L,l,S,s}, an intermediate, letter \v, is also used, which
corresponds to a direction ``reverse'', but this letter does not occur in the output.

The algorithm uses functions:

\centerline{
\texttt{CreateLeft},
\texttt{CreateRight},
\texttt{Reduce},
\texttt{Invert},
\texttt{AlternateCases},
\texttt{CaseForLowerCase},
\texttt{CaseForUpperCase}.
}
These are given as follows.
\begin{enumerate}
\item The function \texttt{CreateLeft}
     produces a left boundary rule with backtracking.  This takes
      $\sigma(\A)$ and replaces $\A,\B,\+,\-$ by $\R,\R,\s,\v$ respectively.  

     \begin{minipage}[t]{11cm}
     \begin{algorithm}[H]
 \DontPrintSemicolon
  \SetKwFunction{FMain}{CreateLeft}
  \SetKwProg{Fn}{Function}{:}{}
  \Fn{\FMain{$ListIn$}}{
    $listL$=[\ ]\;
    \For{$i \leftarrow 1$ \KwTo $n$}{
        \uIf{$x_i=\A$ or $x_i=\B$}{
            $listL.append(\R)$ 
        }
            \uElseIf{$x_i=\+$}{ 
            $listL.append(\s)$
            }
            \uElseIf{$x_i=\-$}{ 
            $listL.append(\v)$
            }            
    }
    \KwRet $ListL$\;
    }
     \end{algorithm}
\end{minipage}
     
     \item The function \texttt{CreateRight}
     produces a right boundary rule with backtracking.  This takes
      $\sigma(\A)$ and replaces $\A,\B,\+,\-$ by $\L,\L,\v,\s$ respectively.

          \begin{minipage}[t]{11cm}
     \begin{algorithm}[H]
 \DontPrintSemicolon
  \SetKwFunction{FMain}{CreateRight}
  \SetKwProg{Fn}{Function}{:}{}
      $listR$=[\ ]\;
  \Fn{\FMain{$ListIn$}}{
    \For{$i \leftarrow 1$ \KwTo $n$}{
        \uIf{$x_i=\A$ or $x_i=\B$}{
            $listR.append(\L)$ 
        }
            \uElseIf{$x_i=\+$}{ 
            $listR.append(\v)$
            }
            \uElseIf{$x_i=\-$}{    
            $listR.append(\s)$
            }
    }
        \KwRet $ListR$\;
    }
     \end{algorithm}
     \end{minipage}

\item
The function $\texttt{Reduce}$, has values as in Table~\ref{tab:red}.
\begin{table}
  $$
  \begin{array}{lllllllllll}
 \texttt{RvR}&\mapsto& \s  &&  \texttt{LvL}&\mapsto& \s&&    \texttt{RvL}&\mapsto& \v\\
  \texttt{svR}&\mapsto& \L  &&  \texttt{svL}&\mapsto&\R&&  \texttt{svs}&\mapsto& \v    \\
  \texttt{Rvs}&\mapsto& \L  &&  \texttt{Lvs}&\mapsto& \R&&  \texttt{LvR}&\mapsto& \v  
  \end{array}
  $$
  \caption{Values of a function \texttt{Reduce} which replaces a three letter subsequence
    with a single letter.}
  \label{tab:red}
\end{table}
\item
  The function
  \texttt{Reduce}
  is extended to apply to sequences of any length, not just length three.

          \begin{minipage}[t]{11cm}
\begin{algorithm}[H]
 \DontPrintSemicolon
  \SetKwFunction{FMain}{Reduce}
  \SetKwProg{Fn}{Function}{:}{}
  \Fn{\FMain{$ListIn$}}{
    $listV = [i \in \text{length}(ListIn): x_i=\v$]\;
   \While{$listV\not=\varnothing$}{
      $j = listV[0]$\;      
Replace the $j-1$ to $j+1$ elements of $ListIn$ 
by $\texttt{Reduce}(ListIn[j-1]ListIn[j]ListIn[j+1])$, the function
as in Table~\ref{tab:red}.\;    
    $listV = [i \in \text{length}(ListIn): x_i=\v$]\;
      }
        \KwRet $ListIn$\;
  }
\end{algorithm}
     \end{minipage}

\item
The function \texttt{Invert} is given by
$$\texttt{Invert}(x_1,x_2,\dots,x_n)=\tilde x_n,\tilde x_{n-1},\dots,\tilde x_1,$$
where $\tilde x$ is an operation which switches $\R\leftrightarrow \L$,
$\r\leftrightarrow \l$,
with case immaterial at this point.

\item
  The function  \texttt{AlternateCases} alternates the case of elements of sequence,
  given a starting case.
  The case of a symbol following $\L,\R,\r,\l$ is the opposite of the case
  of the preceding symbol; the case of a symbol following {\S} or {\s} is the
  same as the case of the preceding symbol.
  This function requires functions \texttt{Case}, \texttt{setCase}
  \texttt{OppositeCase}, which return the case of a letter, set the case of
  a symbol to a given case ,
  and return the opposite case, respectively.
  The case of a letter corresponds to the parity of the square the corresponding
  path starts in, with upper case for even, lower case for odd.

     \begin{minipage}[t]{15cm}
       \begin{algorithm}[H]
             \SetKwInOut{KwIn}{Input}
    \SetKwInOut{KwOut}{Output}
    \KwIn{A list $[x_i]$, $i=1,\dots,n$, with elements in $\{\texttt{L,R,S}\}$, and an initial case,
    \texttt{Upper} or \texttt{Lower}}
    \KwOut{A list $[y_i]$, $i=1,\dots,n$, with $y_i=x_i$ up to case; elements
    are in  $\{\texttt{L,R,S,l,r,s}\}$.}
 \DontPrintSemicolon
  \SetKwFunction{FMain}{AlternateCases}
  \SetKwProg{Fn}{Function}{:}{}
  \Fn{\FMain{$ListIn$,\texttt{InitialCase}}}{
    \texttt{SetCase}($x_1$,\texttt{InitialCase})\;

    \uIf{$x_{i-1}\in\{\texttt{R,L,r,l}\}$
    }{
          \texttt{SetCase}($x_i$,\texttt{OppositeCase}($x_{i-1}$))\;
        }
    \uElseIf{$x_{i-1}\in\{\texttt{S,s}\}$
    }{
         \texttt{SetCase}($x_i$,\texttt{Case}($x_{i-1}$))\;      
            }
  
      \KwRet $ListIn$\;
      }
\end{algorithm}
     \end{minipage}

   \item Funtions
     \texttt{CaseForUpperCase}($\sigma(\A)$) and \texttt{CaseForLowerCase}($\sigma(\A)$)
       to determine the initial case of sequences.

    \begin{minipage}[t]{7cm}
       \begin{algorithm}[H]
             \SetKwInOut{KwIn}{Input}
    \SetKwInOut{KwOut}{Output}
    \KwIn{A list $\sigma(\A)$}
    \KwOut{A case, \texttt{Upper} or \texttt{Lower}}
 \DontPrintSemicolon
  \SetKwFunction{FMain}{CaseForUpperCase}
  \SetKwProg{Fn}{Function}{:}{}
  \Fn{\FMain{$\sigma(\A)$}}{
        \uIf{first term of $\sigma(\A)$ is \A}
  {\KwRet \texttt{Upper}\;}
  \uElse {\KwRet \texttt{Lower}\;}
  }
       \end{algorithm}
       \end{minipage}
           \begin{minipage}[t]{7cm}
       \begin{algorithm}[H]
             \SetKwInOut{KwIn}{Input}
    \SetKwInOut{KwOut}{Output}
    \KwIn{A list $\sigma(\A)$}
    \KwOut{A case, \texttt{Upper} or \texttt{Lower}}
 \DontPrintSemicolon
  \SetKwFunction{FMain}{CaseForLowerCase}
  \SetKwProg{Fn}{Function}{:}{}
  \Fn{\FMain{$\sigma(\A)$}}{
        \uIf{last term of $\sigma(\A)$ is \A}
  {\KwRet \texttt{Lower}\;}
  \uElse {\KwRet \texttt{Upper}\;}
  }
\end{algorithm}
     \end{minipage}

\item
The symbol $+$ is used for the concatenation of lists.

\end{enumerate}
These functions are used in Algorithm~\ref{alg:1},
which results in an L-system defined by the function $\tau$.

\begin{table}
$
  \begin{array}{ll}
    \hline
1& \A\mapsto \sigma(\A)=\A{\myminus}\B\\
 &\B\mapsto \sigma(\B)=\A\+\B\\
\hline
2&\A\mapsto \sigma(\A)= \A{\myminus}\B \+ \A{\myminus}\B \+ \A \+ \B{\myminus}\A \+ \B \+ \A\\
 &\B\mapsto \sigma(\B)= \B{\myminus}\A{\myminus}\B \+ \A{\myminus}\B{\myminus}\A \+ \B{\myminus}\A \+ \B\\
\hline
3 &\A\mapsto \sigma(\A)=\B\+\A{\myminus}\B{\myminus}\A\+\B\+\A\+\B{\myminus}\A\+\B\+\A{\myminus}\B{\myminus}\A{\myminus}\B\+\A{\myminus}\B\+\A\+\B\\
&\B\mapsto \sigma(\B)=\A{\myminus}\B{\myminus}\A \+ \B{\myminus}\A \+ \B \+ \A \+ \B{\myminus}\A{\myminus}\B \+ \A{\myminus}\B{\myminus}\A{\myminus}\B \+ \A \+ \B{\myminus}\A\\
\hline
\hline
\end{array}
$
  \caption{The L-system functions for  three examples of plane-filling folding curves (together with $\+\mapsto\+$ and
    ${\myminus}\mapsto{\myminus}$).}
\label{tab:fill}
\end{table}

\begin{algorithm}
    \SetKwInOut{KwIn}{Input}
    \SetKwInOut{KwOut}{Output}

    \KwIn{A list $[x_i]$, $i=1, 2, \cdots, 2n-1$, corresponding to (and denoted by)  $\sigma(\A)$ where 
    $x_i\in\{\texttt{A,B,+,-}\}$, satisfying (\ref{eqn:1})}
    \KwOut{6 lists, with elements in $\{\texttt{L,R,l,r,S,s}\},$
      corresponding to $\tau(\L),\tau(\l),\tau(\R),\tau(\r),\tau(\S),\tau(\s)$.  A letter
    \v is also used in the algorithm in intermediate step steps.}

    \tcc{Step 1: create a left boundary rule, for $\sigma(\A)$,
      corresponding to the image of a right
      turn, with  backtracking.  This takes
      $\sigma(\A)$ and replaces $\A,\B,\+,\myminus$ by $\R,\R,\s,\v$ respectively. This path starts from an even square.}    
    $listL = \texttt{CreateLeft}(\sigma(\A))$
    
    \tcc{Step 2: removal of backtracking (\v) in left boundary (right turn image),
      by replacing any three element
    subsequence $XvY$ by a single element, $\texttt{Reduce}(XvY)$ as in  Table~\ref{tab:red}}

    $listL=\texttt{Reduce}(listL)$

    \tcc{Step 3: removal of alternate \s terms (assumes index starts at $1$)}

    $listL$ = [$listL$[i] for $i=1,3,\dots,$length(listL) ]\;

    \tcc{Step 4: create a right boundary rule, for $\sigma(\A)$,
      with  backtracking (image of initial left turn).
      This takes
      $\sigma(\A)$ and replaces $\A,\B,\+,\myminus$ by $\R,\R,\s,\v$ respectively. }    
    $listR = \texttt{CreatRight}(\sigma(\A))$

    \tcc{Step 5: removal of backtracking (\v) in right boundary, by replacing three element
    substrings containing $v$ as in Table~\ref{tab:red}}

    $listR=\texttt{Reduce}(listR)$

    \tcc{Step 6: removal of alternate \s terms (assumes index starts at $1$)}

        $listR$ = [$listR$[i] for $i=1,3,\dots,$length(listR) ]\;
    
    \tcc{Step 7: Create left and right boundaries which start from odd squares,
    for the left and right boundaries of $\sigma(\B)$}
    $listr = \texttt{Invert}(listL)$\;
    $listl = \texttt{Invert}(listR)$\;

\tcc{Step 8: Create rule for straight boundary segments (concatenate then reduce).}
    
    $listS = listR + [\v] + listL$\;
    $lists = listL + [\v] + listR$\;
    
    $listS = \texttt{Reduce}(listS)$\;
    $lists = \texttt{Reduce}(lists)$\;

\tcc{Step 9: Adjust the case (upper or lower) of elements of the sequences.
}

\texttt{CaseU} = \texttt{CaseForUpperCase($\sigma(\A)$)};
\texttt{CaseU} = \texttt{CaseForLowerCase($\sigma(\A)$)}

    $listL$=\texttt{AlternateCases}($listl$,\texttt{CaseU});
    $listR$=\texttt{AlternateCases}($listr$,\texttt{CaseU});
    $listS$=\texttt{AlternateCases}($lists$,\texttt{CaseU});
    $listl$=\texttt{AlternateCases}($listl$,\texttt{CaseL});
    $listr$=\texttt{AlternateCases}($listr$,\texttt{CaseL});
    $lists$=\texttt{AlternateCases}($lists$,\texttt{CaseL})
    
    \KwRet{$\tau(\R)=listL,
      \tau(\L)=listR,
      \tau(\S)=listS,
      \tau(\r)=listl,
      \tau(\l)=listr,
      \tau(\s)=lists
      $}
    \tcc{note that listL, which gives a description of the left side of $\sigma(\A)$,
      is the word for $\tau(\R)$, since the left side of {\A} is given by a single right turn,
      as illustrated in Figure~\ref{fig:3examples}.  Similarly for the other terms.}
    \caption{Computation of the L-system for the boundary of a plane-filling
      folding curve}
    \label{alg:1}
\end{algorithm}

\section{Examples}
Before giving a proof that the algorithm does indeed construct an L-system for the
boundary of the given plane-filling curve,
we show how this algorithm applies in three examples.
We do not discuss how to find examples of folding curve L-systems, which
can be found in \cite{arndt} and \cite{AH}.  Case (1) is the Heighway dragon and case (2) is from
Figure 1.4.1-D from \cite{arndt}.
In Tables~\ref{tab:fill} and \ref{tab:bou}, the L-systems, $\sigma$ and the corresponding $\tau$ are shown, along with the
curves corresponding to
$\sigma(\A)$, $\sigma^3(\A)$ and $\tau(\R), \tau(\L), \tau^3(\R), \tau^3(\L)$
(possibly drawn at different scales and orientations to each other).

\begin{table}
$$
\begin{array}{llr||l|l|l|l}
 && & 1 & 2 & 3\\
  \hline
  \L&\mapsto&\tau(\L)= & \L\l   & \texttt{LsrR}     & \texttt{rLsr} \\
  \R&\mapsto&\tau(\R)= & \S  & \texttt{SSRsrR}   & \texttt{rLrRslRr}\\
  \l&\mapsto&\tau(\l)= & \S  & \texttt{lLslSS}   & \texttt{LlRslLrL}\\
  \r&\mapsto&\tau(\r)= & \R\r   & \texttt{lLsr}     & \texttt{LsrL}\\
  \S&\mapsto&\tau(\S)= & \L\r  & \texttt{SSS}      & \texttt{rLrRslLrL}\\
  \s&\mapsto&\tau(\s)= & \R\l  & \texttt{lLslSRsrR}& \texttt{LlRslRr}\\
  \end{array}
$$
\caption{
  The L-systems for the boundaries of the curves in Table~\ref{tab:fill}.}
\label{tab:bou}
\end{table}

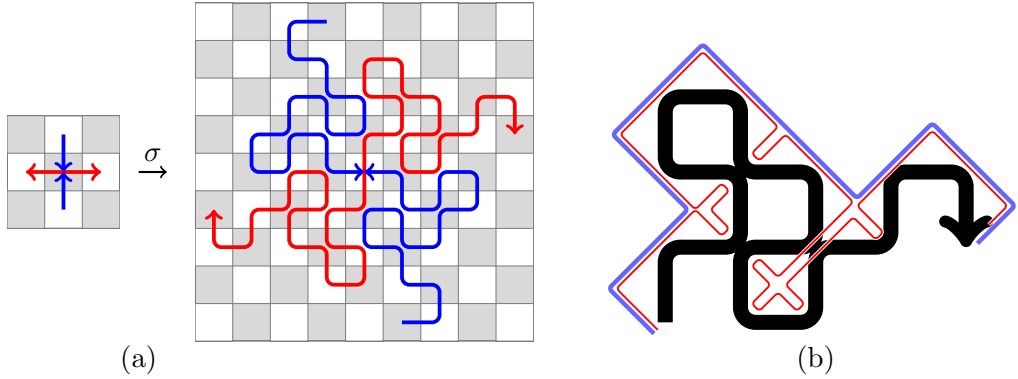
\begin{figure}
  \centering
  \begin{tikzpicture}
    \begin{scope}[scale={0.5},shift={(4,4)}]

      \begin{scope}[shift={(-8,0)}]
        \clip (-1.5, -1.5) rectangle (1.5,1.5);
        \draw[fill=white](-1.5, -1.5) rectangle (1.5,1.5);
        \foreach\i in {-1,...,1}{
        \foreach\j in {-1,...,1}{
          \draw[gray,fill=gray!30!white](\i+\j-0.5,\i-\j-0.5) rectangle (\i+\j+0.5,\i-\j+0.5);
          }}
\draw[red,line width = 0.5mm,rounded corners=1mm,->](0,0)--(1,0);
\begin{scope}[rotate={90}]\draw[blue,line width = 0.5mm,rounded corners=1mm,<-](0,0)--(1,0);\end{scope}
  \begin{scope}[rotate={180}]\draw[red,line width = 0.5mm,rounded corners=1mm,->](0,0)--(1,0);\end{scope}
\begin{scope}[rotate={270}]\draw[blue,line width = 0.5mm,rounded corners=1mm,<-](0,0)--(1,0);\end{scope}          
\end{scope}

\draw[thick,->](-6,0)--node[above]{$\sigma$}(-5.25,0);

\begin{scope}
  \clip (-4.5, -4.5) rectangle (4.5,4.5);
        \draw[fill=gray!30!white](-4.5, -4.5) rectangle (4.5,4.5);
        \foreach\i in {-4,...,4}{
          \foreach\j in {-4,...,4}{
                      \draw[gray,fill=white](\i+\j-0.5,\i-\j-0.5) rectangle (\i+\j+0.5,\i-\j+0.5);
          }}
\end{scope}
\draw[red,line width = 0.5mm,rounded corners=1mm,->]\mypath
\begin{scope}[rotate={90}]\draw[blue,line width = 0.5mm,rounded corners=1mm,<-]\mypath\end{scope}
  \begin{scope}[rotate={180}]\draw[red,line width = 0.5mm,rounded corners=1mm,->]\mypath\end{scope}
\begin{scope}[rotate={270}]\draw[blue,line width = 0.5mm,rounded corners=1mm,<-]\mypath\end{scope}
    \end{scope}

\begin{scope}[shift={(6,0)}]
  \draw[black,line width = 2mm,rounded corners=2mm,->]\mypath
      \draw[rounded corners=0.5mm,white, line width = 0.5mm](-0.1,-0.1)--++(-0.55,0.55)--++(1.05,1.05)
       --++(0.35,-0.35)--++(0.1,0.1)--++(-0.25,0.25)--++(0.25,0.25)--++(-0.1,0.1)
       --++(-0.35,-0.35)--++(-1.0,1.0)
       --++(1.1,1.1)--++(1.,-1.)
       --++(-0.35,-0.35)--++(0.1,-0.1)--++(0.35,0.35)
       --++(.9,-.9)--++(-1,-1)
       --++(-0.25,0.25)--++(-0.1,-0.1)--++(0.25,-0.25)
       --++(-0.25,-0.25)--++(0.1,-0.1)--++(0.25,0.25)
       --++(0.25,-0.25)--++(0.1,0.1)--++(-0.25,0.25)
       --++(0.9,0.9)--++(0.3,-0.3)--++(0.1,0.1)--++(-0.35,0.35)
       --++(1,1)--++(1,-1)--++(-0.25,-0.25);

       \draw[rounded corners=0.5mm, line width=0.25mm,red](-0.1,-0.1)--++(-0.55,0.55)--++(1.05,1.05)
       --++(0.35,-0.35)--++(0.1,0.1)--++(-0.25,0.25)--++(0.25,0.25)--++(-0.1,0.1)
       --++(-0.35,-0.35)--++(-1.0,1.0)
       --++(1.1,1.1)--++(1.,-1.)
       --++(-0.35,-0.35)--++(0.1,-0.1)--++(0.35,0.35)
       --++(.9,-.9)--++(-1,-1)
       --++(-0.25,0.25)--++(-0.1,-0.1)--++(0.25,-0.25)
       --++(-0.25,-0.25)--++(0.1,-0.1)--++(0.25,0.25)
       --++(0.25,-0.25)--++(0.1,0.1)--++(-0.25,0.25)
       --++(0.9,0.9)--++(0.3,-0.3)--++(0.1,0.1)--++(-0.35,0.35)
       --++(1,1)--++(1,-1)--++(-0.25,-0.25);

       \draw[rounded corners=0.5mm, line width=0.6mm,blue!60!white](-0.15,-0.15)--++(-0.6,0.6)
       --++(1.05,1.05) --++(-1.0,1.0)
       --++(1.2,1.2)--++(2.05,-2.05) --++(1.,1.) --++(1.1,-1.1)
       --++(-0.5,-0.5);
       
\end{scope}
\node at (-1,-.5){(a)};
\node at (8,-0.5){(b)};
    \end{tikzpicture}
  \caption{
Example 3, with {\tiny $\sigma(\A)=\texttt{\B+A-\B-A+\B+A+\B-A+\B+A-\B-A-\B+A-\B+A+\B}$}, from Table~\ref{fig:3examples}.
(a) Application of L-system operation once to four line segments, two {\A} and two {\B} direction.
Gray squares are even.
(b)    Diagram of left side boundary, in blue.  Red path is prior to removal of backtracking.
\texttt{A, R, L} start in odd squares. 
  }
\label{fig:part3A+1}
  \end{figure}

Starting with the L-system given by example 3 of Table~\ref{tab:fill},
Figure~\ref{fig:part3A+1}, right shows: the path $\sigma(\A)$, in black;
the path corresponding to \texttt{CreateRight}$(\sigma(\A))$, in red;
the reduced path, obtained from Algorithm~\ref{alg:1}, Steps 3,4, in blue.
Application of Steps 4, 5, 6, 9 of Algorithm~\ref{alg:1} are as in Table~\ref{tab:exam}
\begin{table}
$
\arraycolsep=1.4pt
\begin{array}{l|lcccccccccccccccccccccccccccccccccc}
0&  \text{ original curve} &  \B&\+&\A&\myminus&\B&\myminus&\A&\+&\B&\+&\A&\+&\B&\myminus&\A&\+&\B&\+&\A&\myminus&\B&\myminus&\A&\myminus&\B&\+&\A&\myminus&\B&\+&\A&\+&\B\\
    \hline
4&  \text{ right boundary} &  \L&\v&\L&\s&\L&\s&\L&\v&\L&\v&\L&\v&\L&\s&\L&\v&\L&\v&\L&\s&\L&\s&\L&\s&\L&\v&\L&\s&\L&\v&\L&\v&\L\\
    \hline
5  &&   &\s& &\s&\L&\s& &\s& &\v& &\s& &\s& &\s& &\v&\L&\s&\L&\s&\L&\s& &\s& &\s& &\s& &\v&\L\\
&  &   &\s& &\s&\L&\s& & & &\v& & & &\s& & & &\R& &\s&\L&\s&\L&\s& &\s& &\s& & & &\R& \\
&  &   &\s& &\s&\L& & & & &\v& & & & & & & &\R& &\s&\L&\s&\L&\s& &\s& &\s& & & &\R& \\
&  &   &\s& &\s& & & & & &\v& & & & & & & & & &\s&\L&\s&\L&\s& &\s& &\s& & & &\R& \\
&  &   &\s& & & & & & & &\v& & & & & & & & & & &\L&\s&\L&\s& &\s& &\s& & & &\R& \\
&    &   & & & & & & & & &\R& & & & & & & & & & & &\s&\L&\s& &\s& &\s& & & &\R& \\
\hline
6&        &   & & & & & & & & &\R& & & & & & & & & & & & &\L& & &\s& & & & & &\R& \\
    \hline
9&            &   & & & & & & & & &\r& & & & & & & & & & & & &\L& & &\s& & & & & &\r& \\
        \hline
\end{array}
$
\caption{Application of Algorithm~\ref{alg:1} steps 4,5,6,9 to obtain $\tau(\L)$.}
  \label{tab:exam}
\end{table}
Thus in this example we obtain
$$\tau:\L\mapsto  \texttt{rLsr}.$$
Note that in this example, $\sigma(\A)$ starts with {\B}, so since {\A} starts in
even squares and {\B} in odd squares, the central square in Figure~\ref{fig:part3A+1}
changes from even to odd.

\section{Proof that Algorithm~\ref{alg:1} produces the boundary
  of the plane-filling curve}

To prove that
Algorithm~\ref{alg:1} has the required output,
we justify that:

\begin{itemize}
\item $\texttt{CreateLeft}(\sigma(\A))$
  and $\texttt{CreateRight}(\sigma(\A))$
  create left and right boundaries of $\sigma(\A)$,
  with backtracking.
\item
  The reduction function creates a boundary with no backtracking.
\item
  The paths $\tau(\S)$ and $\tau(\s)$ are also
  boundary components.
\item
  The rule for the alteration of the cases of letters in
  the sequences $\tau$ is correct.
\end{itemize}
From these, we can see that
$\tau^n(\R)$ union $\tau^n(\L)$
form the boundary of
$\sigma^n(\R)$. The final result, that the limits of
$\tau^\infty(\R)$ union $\tau^\infty(\L)$
form the boundary of
$\sigma^\infty(\R)$ follows by induction.
Similarly for the letters $\B, \r, \l$.

\subsection{Creation of boundary paths}

\begin{figure}
\scalebox{0.8}{
  \begin{tikzpicture}

    \node at (-0.2,2.7){$\tau_0(\R)$};
    \node at (-0.2,1.3){$\tau_0(\L)$};
        \node at (-0.65,2){$\sigma_0(\A)$};

    \draw[line width = 3mm,->](0,2)--(9,2);
    \draw[red,rounded corners = 1mm,line width = 2mm,->]
    (0,2.2)
    \foreach\i in {0,...,3}{
--++(1,1)--++(1,-1)
    }
    ;
    \draw[blue,rounded corners = 1mm,line width = 2mm,->]
    (0,1.8)
    \foreach\i in {0,...,3}{
--++(1,-1)--++(1,1)
    }
    ;

    \foreach\i in {0,...,4}{
\draw[fill=white] (2*\i,2) circle (0.2);
\node at (2*\i,2){?};
    }

   \foreach\i in {0,...,3}{
\node at (2*\i+1,3.4){\R};
\node at (2*\i+1,.6){\L};
   }

   \node at (4,0){(a)};
   \node at (10,0){(b)};
   \node at (14,0){(c)};
   
\begin{scope}[shift={(10,0)}]
      
  \draw[line width = 3mm,->,rounded corners = 1mm](0,2)--(2,2)--(2,0);
  \draw[red,rounded corners = 1mm,line width = 2mm,->]
  (0,2.3)--++(1.05,1.05)--++(2.3,-2.3)--(2.3,-0.2);
\draw[fill=white] (1.9,1.9) circle (0.2);
\node at(1.9,1.9){\R};

\draw[blue,rounded corners = 1mm,line width = 2mm,->]
(0,1.7)--++(0.8,-0.8)--++(0.8,0.8)--++(0.2,-0.2)
--++(-0.8,-0.8)--++(0.8,-0.8);

\draw[fill=white] (2.2,2.2) circle (0.2);
\node at(2.2,2.2){\S};

\draw[fill=white] (1.6,1.6) circle (0.2);
\node at(1.6,1.6){v};

\node at(1,3.5){\R};
\node at(3,1){\R};

\node at(0.7,0.6){\L};
\node at(0.9,0.4){\L};

    \node at (-0.2,2.7){\tiny $\tau_1(\R)$};
    \node at (-0.2,1.3){\tiny $\tau_1(\L)$};
        \node at (-0.4,2){\tiny $\sigma(\A)$};

\end{scope}

\begin{scope}[shift={(14,3)},yscale = {-1}]

\node at(1,3.5){\L};
\node at(3,1){\L};

\node at(0.7,0.6){\R};
\node at(0.9,0.4){\R};

  \draw[line width = 3mm,->,rounded corners = 1mm](0,2)--(2,2)--(2,0);
  \draw[blue,rounded corners = 1mm,line width = 2mm,->]
  (0,2.3)--++(1.05,1.05)--++(2.3,-2.3)--(2.3,-0.2);
\draw[fill=white] (1.9,1.9) circle (0.2);
\node at(1.9,1.9){\L};

\draw[red,rounded corners = 1mm,line width = 2mm,->]
(0,1.7)--++(0.8,-0.8)--++(0.8,0.8)--++(0.2,-0.2)
--++(-0.8,-0.8)--++(0.8,-0.8);

\draw[fill=white] (2.2,2.2) circle (0.2);
\node at(2.2,2.2){\S};

\draw[fill=white] (1.6,1.6) circle (0.2);
\node at(1.6,1.6){\v};

  \end{scope}

  \end{tikzpicture}
  }
  \caption{Justification for the description of
    $\texttt{CreateLeft}(\A)$
    and $\texttt{CreateRight}(\A)$.
    In (a) a straight black line represents the to
    be folded plane-filling curve, where ``?'' are to be determined,
    either $\+$ (right) or $\-$ (left).
    The red curve $\tau_0(R)$
    above the black curve, is the left side of the black curve.
        The blue curve $\tau_0(\L)$
        below the black curve, is the right side of the black curve.
        These curves consist of segments with alternating turns to be determined,
        depending on the black curve.
        (b) shows the effect of a right turn of the folding curve on the boundary curves,
        and (c) shows the effect of a left turn.
  }
  \label{fig:bends}
\end{figure}
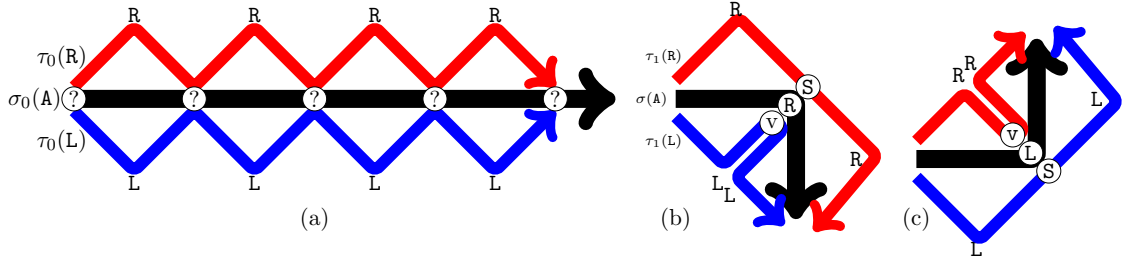

The \texttt{CreateLeft} and \texttt{CreateRight} functions
correspond geometrically to placing a right angle triangle
along each edge of the original curve, $\sigma(\A)$, on either the
left or right side.
For brevity, we write
$\tau_1(\R)=\texttt{CreateLeft}(\A)$ and
$\tau_1(\L)=\texttt{CreateRight}(\A)$.
Let
$\tau_0(\R)$,
$\tau_0(\L)$, and
$\sigma_0(\A)$ be paths shown in
Figure~\ref{fig:bends} (a), corresponding to $\sigma(\A)$ but with all turns
(letters {\+} and {\myminus}) removed.
The result of inserting {\+} or {\-} has the effect shown
in Figure~\ref{fig:bends} (b) and (c) at any particular junction.
Considering $\tau_0(\R)$ to the left of $\sigma_0(\A)$,
and  $\tau_0(\L)$
 to the right of $\sigma_0(\A)$, before applying the
folding directions, so that curves do not touch, except at mid points of squares, and they never cross.
Given that the turns are never more than $90$, and that given that
$\sigma(\A)$ is assumed to describe a curve that is non self intersecting,
then we must have that after applying the folding sequence,
$\tau_1(\R)$ is still to the left, and $\tau_1(\L)$ is still to the right.
In particular, this can also be seen to entail that
the right side of the path
$\tau_1(\R)$ always touches the curve $\sigma(\A)$; the left side of
$\tau_1(\R)$ either touches $\sigma(\A_1)$ or $\sigma(\B_1)$ where $\A_1$ and
$\B_1$ are different iterations of the plane-filling curve, starting from a different edge
that the particular $\A$ which $\R$ is left of, or else,
segments of the left side of $\tau_1(\R)$ touch the curve $\tau_1(\R)$.
Similarly for the left side boundary.
A right turn of $\sigma(\A)$ corresponds to a straight ahead instruction for $\tau_1(\R)$,
and a reverse instruction for $\tau_1(\L)$, and conversely for a left turn.
The explains the construction of the initial \texttt{CreateLeft} and \texttt{CreateRight} functions.
Here the notion of ``touch'' is interpreted to mean touching at the mid-points of squares.
Provided curves are drawn from centre to centre, they will intersect at the centre
points of squares.
In the limit, when the squares become arbitrarily small,
the limits of the curves will touch continuously.

\subsection{Reduction of paths, removing backtracking}

\begin{figure}
  \scalebox{0.7}{
    \begin{tikzpicture}
      \filldraw[gray!10!white](0,0) rectangle (2,2);
    \draw[line width = 1mm,->,rounded corners = 1mm] (0.1,0)--++(0.9,0.9)--++(0.9,-0.9)--++(.1,0.1)--++(-0.9,0.9)--++(0.9,0.9);
    \node at(0.8,1){\R};
    \node at(1.1,1.3){\R};
    \node at(2.1,0){\v};
    \node at (1,-0.4){$\texttt{RvR}\mapsto \S$};
    \begin{scope}[shift={(2.5,0)}]
                      \filldraw[gray!10!white](0,0) rectangle (2,2);
    \draw[line width = 1mm,->,rounded corners = 1mm] (0.1,0)--++(0.9,0.9)--++(0.9,-0.9)--++(.1,0.1)--++(-0.9,0.9)--++(-0.9,0.9);
    \node at(0.8,1){\R};
    \node at(1.1,1.3){\S};
    \node at(2.1,0){\v};
    \node at (1,-0.4){$\texttt{RvS}\mapsto \L$};
\end{scope}
    \begin{scope}[shift={(5,0)}]
                      \filldraw[gray!10!white](0,0) rectangle (2,2);
    \draw[line width = 1mm,->,rounded corners = 1mm] (0.1,0)--++(0.9,0.9)--++(0.9,-0.9)--++(.1,0.1)--++(-1,1)--++(-1,-1);
    \node at(1,0.5){\R};
    \node at(1.1,1.4){\L};
    \node at(2.1,0){\v};
    \node at (1,-0.4){$\texttt{RvL}\mapsto \v$};
\end{scope}
    \begin{scope}[shift={(7.5,0)}]
                      \filldraw[gray!10!white](0,0) rectangle (2,2);
    \draw[line width = 1mm,->,rounded corners = 1mm] (0,0.1)--++(0.9,0.9)--++(-0.9,0.9)--++(.1,0.1)--++(0.9,-0.9)--++(0.9,0.9);
    \node at(1,0.8){\L};
    \node at(1.3,1.1){\L};
    \node at(-0.2,2){\v};
    \node at (1,-0.4){$\texttt{LvL}\mapsto \S$};
    \end{scope}
    \begin{scope}[shift={(10,0)}]
                      \filldraw[gray!10!white](0,0) rectangle (2,2);
    \draw[line width = 1mm,->,rounded corners = 1mm] (0,0.1)--++(0.9,0.9)--++(-0.9,0.9)--++(.1,0.1)--++(0.9,-0.9)--++(0.9,-0.9);
    \node at(1,0.8){\L};
    \node at(1.1,1.3){\S};
        \node at(-0.2,2){\v};
    \node at (1,-0.4){$\texttt{RvS}\mapsto \R$};
\end{scope}
    \begin{scope}[shift={(12.5,0)}]
                      \filldraw[gray!10!white](0,0) rectangle (2,2);
    \draw[line width = 1mm,->,rounded corners = 1mm] (0,0.1)--++(0.9,0.9)--++(-0.9,0.9)--++(.1,0.1)--++(1,-1)--++(-1,-1);
    \node at(0.5,1){\L};
    \node at(1.3,1){\R};
            \node at(-0.2,2){\v};
    \node at (1,-0.4){$\texttt{LvR}\mapsto \v$};
\end{scope}
    \begin{scope}[shift={(15,0)}]
                      \filldraw[gray!10!white](0,0) rectangle (2,2);
    \draw[line width = 1mm,->,rounded corners = 1mm] (0.1,0)--++(1.9,1.9)--++(-.1,0.1)--++(-0.9,-0.9)--++(-0.9,0.9);
    \node at(1,1.4){\L};
    \node at(1.2,0.8){\S};
    \node at(2,2.1){\v};
    \node at (1,-0.4){$\texttt{SvR}\mapsto \L$};
    \end{scope}
    \begin{scope}[shift={(17.5,0)}]
                      \filldraw[gray!10!white](0,0) rectangle (2,2);
    \draw[line width = 1mm,->,rounded corners = 1mm] (0,0.1)--++(1.9,1.9)--++(.1,-0.1)--++(-0.9,-0.9)--++(0.9,-0.9);
    \node at(1,0.8){\L};
    \node at(1,1.4){\S};
        \node at(2,2.1){\v};
    \node at (1,-0.4){$\texttt{RvS}\mapsto \R$};
\end{scope}
    \begin{scope}[shift={(20,0)}]
                      \filldraw[gray!10!white](0,0) rectangle (2,2);
    \draw[line width = 1mm,->,rounded corners = 1mm] (0,0.1)--++(1.9,1.9)--++(.1,-0.1)--++(-2,-2);
    \node at(0.8,1.1){\S};
    \node at(1.1,0.8){\S};
            \node at(2,2.1){\v};
    \node at (1,-0.4){$\texttt{SvS}\mapsto \v$};
    \end{scope}

\begin{scope}[shift={(0,-3)}]
      \filldraw[gray!10!white](0,0) rectangle (2,2);
    \draw[line width = 1mm,->,rounded corners = 1mm] (0.1,0)--++(1.9,1.9);
\filldraw (1,0.9) circle (0.1);
    \node at(0.8,1){\S};
    \begin{scope}[shift={(2.5,0)}]
                \filldraw[gray!10!white](0,0) rectangle (2,2);
    \draw[line width = 1mm,->,rounded corners = 1mm] (0.1,0)--++(1,1)--++(-0.9,0.9);
    \node at(1.1,1.3){\L};
\end{scope}
\begin{scope}[shift={(5,0)}]
          \filldraw[gray!10!white](0,0) rectangle (2,2);
  \draw[->,line width = 1mm,rounded corners = 1mm]
  (0.1,0)--++(1,1)--++(-0.1,0.1)--++(-1.2,-1.2);
      \node at(1.2,1.2){\v};
\end{scope}
\begin{scope}[shift={(7.5,0)}]
        \filldraw[gray!10!white](0,0) rectangle (2,2);
    \draw[line width = 1mm,->,rounded corners = 1mm] (0.1,0)--++(1.9,1.9);
\filldraw (1,0.9) circle (0.1);
    \node at(0.8,1){\S};
    \end{scope}
\begin{scope}[shift={(10,0)}]
                \filldraw[gray!10!white](0,0) rectangle (2,2);
    \draw[line width = 1mm,->,rounded corners = 1mm] (0,0.1)--++(0.9,0.9)--++(0.9,-0.9);
    \node at(1.3,1.1){\R};
\end{scope}
\begin{scope}[shift={(12.5,0)}]
  \filldraw[gray!10!white](0,0) rectangle (2,2);
  \draw[->,line width = 1mm,rounded corners = 1mm]
  (0.1,0)--++(1,1)--++(-0.1,0.1)--++(-1.2,-1.2);
      \node at(1.2,1.2){\v};
\end{scope}
\begin{scope}[shift={(15,0)}]
                \filldraw[gray!10!white](0,0) rectangle (2,2);
    \draw[line width = 1mm,->,rounded corners = 1mm] (0.1,0)--++(1,1)--++(-0.9,0.9);
    \node at(1.1,1.3){\L};
    \end{scope}
\begin{scope}[shift={(17.5,0)}]
                  \filldraw[gray!10!white](0,0) rectangle (2,2);
    \draw[line width = 1mm,->,rounded corners = 1mm] (0,0.1)--++(0.9,0.9)--++(0.9,-0.9);
    \node at(1.3,1.1){\R};
\end{scope}
\begin{scope}[shift={(20,0)}]
                  \filldraw[gray!10!white](0,0) rectangle (2,2);
  \draw[->,line width = 1mm,rounded corners = 1mm]
  (0.1,0)--++(1,1)--++(-0.1,0.1)--++(-1.2,-1.2);
      \node at(1.2,1.2){\v};
    \end{scope}

\end{scope}

  \end{tikzpicture}
  }
  \caption{Illustration of the backtracking reduction operation.  Top row shows
    the initial configuration, bottom row the resulting path after reduction.  Start
  and end points are unchanged.}
  \label{fig:reducecases}
  \end{figure}
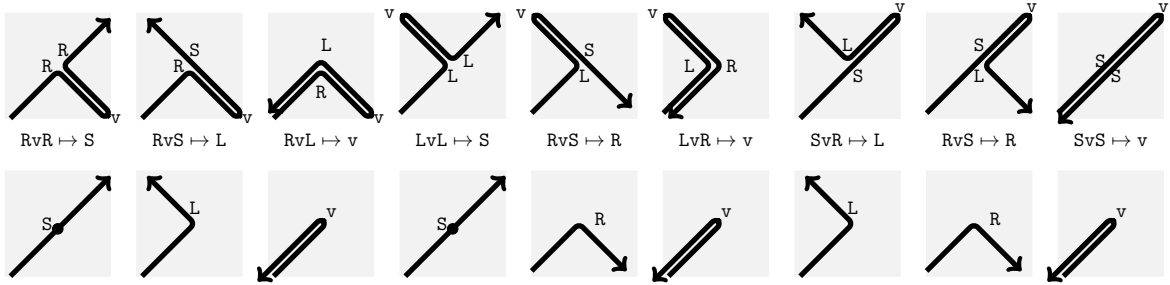

The reduction function corresponds to removal of backtracking.
Figure~\ref{fig:part3A+1} (b), illustrates example (3) from Tables~\ref{tab:bou} and \ref{tab:fill}.
In general, the reduction rules, shown in Table~\ref{tab:red} are justified by the illustrations in
Figure~\ref{fig:reducecases}.
Application of the reduction function results in a path with no reversing.
This is
because each application of $\texttt{Reduce}$ replaces a segment of
the form
$XvZ$ with either $\texttt{L,R,S,v}$.
The only way a {\v} will remain after replacements is if {\v} is an initial or final
term,
or if there are two adjacent {\v} terms in the sequence at some stage.
These cases can not occur. 

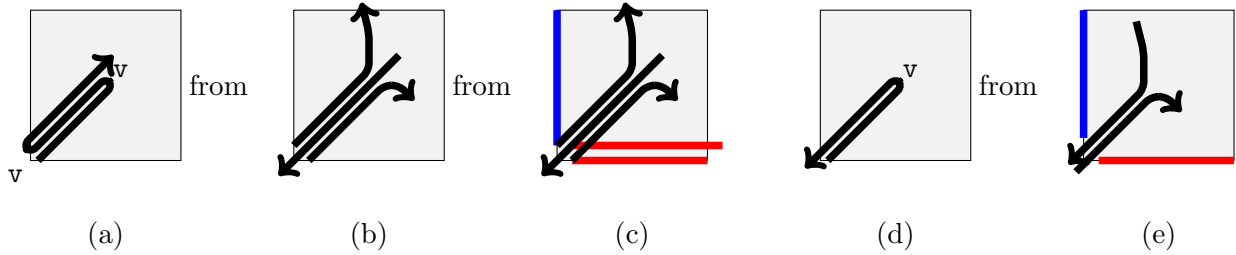
\begin{figure}

  \begin{tikzpicture}
    \draw[fill=gray!10!white] (0,0) rectangle (2,2);
      \draw[->,line width = 1mm,rounded corners = 1mm]
  (0.1,0)--++(1,1)--++(-0.1,0.1)--++(-1,-1)--++(-0.1,0.1)--++(1.2,1.2);
      \node at(1.2,1.2){\v};
            \node at(-.2,-.2){\v};
                        \node at (2.5,1){from};

    \begin{scope}[shift={(3.5,0)}]    
      \draw[fill=gray!10!white] (0,0) rectangle (2,2);
      \draw[->,line width = 1mm,rounded corners = 1mm](0,0.2)--++(1,1)--++(0,0.5)--++(-0.1,.4);
      \draw[<-,line width = 1mm,rounded corners = 1mm](-0.2,-0.2)--++(1.6,1.6);
            \draw[->,line width = 1mm,rounded corners = 1mm](0.2,0)--++(1,1)--++(0.2,0)--++(0.2,-0.2);
            \node at (2.5,1){from};
    \end{scope}

    \begin{scope}[shift={(7,0)}]    
      \draw[fill=gray!10!white] (0,0) rectangle (2,2);
                  \draw[red,line width = 1mm,rounded corners = 1mm](0.2,0)--(2,0);
                        \draw[red,line width = 1mm,rounded corners = 1mm](0.2,0.2)--++(2,0);
                        \draw[blue,line width = 1mm,rounded corners = 1mm](0,0.2)--(0,2);
            \draw[->,line width = 1mm,rounded corners = 1mm](0,0.2)--++(1,1)--++(0,0.5)--++(-0.1,.4);
      \draw[<-,line width = 1mm,rounded corners = 1mm](-0.2,-0.2)--++(1.6,1.6);
            \draw[->,line width = 1mm,rounded corners = 1mm](0.2,0)--++(1,1)--++(0.2,0)--++(0.2,-0.2);

   \begin{scope}[shift={(3.5,0)}]    
    \draw[fill=gray!10!white] (0,0) rectangle (2,2);
      \draw[->,line width = 1mm,rounded corners = 1mm]
  (0.1,0)--++(1,1)--++(-0.1,0.1)--++(-1.2,-1.2);
      \node at(1.2,1.2){\v};
                        \node at (2.5,1){from};
   \end{scope}

   \begin{scope}[shift={(7,0)}]
      \draw[fill=gray!10!white] (0,0) rectangle (2,2);
                  \draw[red,line width = 1mm,rounded corners = 1mm](0.2,0)--(2,0);
                        \draw[blue,line width = 1mm,rounded corners = 1mm](0,0.3)--(0,2);
            \draw[<-,line width = 1mm,rounded corners = 1mm](-0.2,-0.05)--++(1,1)--++(0,0.5)--++(-0.1,.4);
            \draw[->,line width = 1mm,rounded corners = 1mm](-0.1,-0.15)--++(1,1)--++(0.2,0)--++(0.2,-0.2);
      \end{scope}

   \node at (-6,-1){(a)};
   \node at (-2.5,-1){(b)};
   \node at (1,-1){(c)};
   \node at (4.5,-1){(d)};
         \node at (8,-1){(e)};

               \end{scope}

    \end{tikzpicture}
  
  \caption{Reason we can not have $\texttt{vv}$ in the reduced boundary path, and the reduced path can not begin or end with
    $\v$.}
  \label{fig:redvv}
  \end{figure}

We can not have {\v} followed by {\v} in a reduced system.
This would correspond to three paths as
parts of $\texttt{CreateLeft}(\A)$ down the same diagonal, as in Figure~\ref{fig:redvv}.
This is because the reduction step only removes edges, it can not add any, so a situation as in
Figure~\ref{fig:redvv} left must come from some original configuration as in
Figure~\ref{fig:redvv} centre, but each of these diagonal components is associated with
a component of $\sigma(\A)$, by construction,
which must lie along the edges of this square which are adjacent to the corner {\v},
as in Figure~\ref{fig:redvv} right.  But
this would mean $2$ edges of $\sigma(\A)$ along the same space, but
because $\sigma(\A)$ is non self intersecting (by assumption), this is impossible.
Hence we can never reduce to the subsequence
$\v\v$, of $\v$ directly following $\v$, being contained in the reduction of the
left or right boundaries of $\sigma(\A)$.
This means that if there is a {\v} in the sequence
obtained after a reduction step, we can apply another reduction step.
Since the only elements of the sequence are $\texttt{R, L, S, v}$, and
there are only a finite number of terms, eventually, it must be
that we can not apply any more reduction steps,
and at this point, we must have no remaining {\v} terms.

We can not have that the reduced path starts or ends in $\v$.
This is because if this was the case, we would have a sitution as in
Figure~\ref{fig:redvv} (d), coming from
Figure~\ref{fig:redvv} (e), where the image of $\sigma(A)$ contains components $A$ and $B$ both with
start or end points in the initial vertex of $A$.  But this contradicts $\sigma(A)$ not overlapping any other
$\sigma(B)$ with the same start point.  This is illustrated in Figure~\ref{fig:part3A+1} (a), where
the central point is the starting or end point of exactly two copies of $\sigma(A)$ and two of $\sigma(B)$.

Steps $3$ and $6$ remove the excess letters {\s}.
This is because for the L-system $\tau$, the ``forward'' instruction is assumed between
each directional instruction.  The reduction operation removes three terms at a time, of
the form $XyZ$, with $\{X,Z\}\in\{\texttt{L,R,s,v}\}$, and $y\in\{\texttt{v,s}\}$ and replaces
with a symbol $W\in\{\texttt{L,R,s,v}\}$.
Assume we index from $1$, so the sequence $\sigma(\A)$ has the form $x_1,x_2,\dots$.
Then this can be considered an operation which replaces
a term at an odd location with another term at an odd location, and changes the index
of remaining terms, but not their parity, so elements with even index still are either
{\s} or {\v} at any stage.
Given the above argument, after the end of the reduction process, every other term is {\s},
since all letters {\v}  have been eliminated.  Thus every other letter corresponds to a ``move forwards''
instruction, but these are not required in our L-system; they are an assumed instruction
following every directional instruction.
Therefore they can be removed.  (An alternative would be to include another  (fixed) symbol in the
L-system for forward; note that {\S} and {\s} in our system do not mean ``forward'', they can be
considered place holders for the ``don't turn'' instruction.)

Once we have a path with no reversing, obtained from $\tau_1(\R)$ or $\tau_1(\L)$,
we need to justify that this path does not touch itself.
This is because the only way a path could touch itself, but not have a reverse step, would involve
an isolated region of the plane, which does not contain any squares containing a segment of
any $\sigma(\A)$, using the fact that the image of the plane-filling curve is connected.
This is illustrated in Figure~\ref{fig:hole}.
But this contradicts that the iterates $\sigma^n(\A)$ contain components in all
squares of $\mathbb Z^n$, at any stage of iteration, since they are assumed to converge
to plane-filling curves, and this is part of the construction in for example \cite{dekking}.
\begin{figure}
  \scalebox{0.7}{
  \begin{tikzpicture}
    \draw[gray!50!white] (0,0) grid (11,6);
    \draw[red,line width = 2mm, rounded corners=1mm](0,0)--(3.9,3.9)--++(0.99,-0.95)
    --++(-0.9,-0.9)--++(2,-2)--++(1,1)--++(1,-1)--++(2,2)--++(-2,2)--++(-1.9,-1.9)--++(-3,3)--++(1,1);

    \draw[gray,line width = 1mm, rounded corners=2mm]((0,1.5)--++(0.5,0)--++(0,1)
    --++(1,0)--++(0,-1)--++(-1,0)--++(0,-1)--++(-0.5,0);

        \draw[gray,line width = 1mm, rounded corners=2mm]((0.5,6)--++(0,-0.5)
        --++(1,0)--++(0,-1)--++(1,0)--++(0,-1)--++(1,0)--++(0,1)--++(-1,0)--++(0,1)--++(1,0)--++(0,0.5);

                \draw[black,line width = 1mm, rounded corners=3mm]((9.5,6)--++(0,-0.5)
                --++(-1,0)--++(0,-1)--++(1,0)--++(0,-1)--++(1,0)--++(0,-1)--++(-1,0)--++(0,1)--++(-1,0)--++(0,1)
                --++(-1,0)--++(0,-1)--++(-1,0)--++(0,-1)--++(-1,0)--++(0,1)--++(-1,0)--++(0,0.5);

                \node at(7.5,1.5){\parbox{5cm}{\small Empty region that can not contain any
                  $\sigma(\A)$ for {\A} left or right of \R, which is impossible}};

\node at (5,5){$\tau(\R)$};
              \node at(4,0.5){\parbox{5cm}{segment where $\tau(\R)$\\ touches itself}};
              \draw[->](2,0.8)--++(2.4,2.4);
              \draw[->](4.5,5)--(3.4,5);
                \node at(10,6){part of $\sigma(\A)$ for {\A} to the right
                  of $\tau(\R)$};
                                \node at(0.5,4){\parbox{4cm}{part of $\sigma(\A)$ for some \\{\A} to the left
                  of $\tau(\L)$}}; 
    
  \end{tikzpicture}
  }
\caption{Example situation if there is no backtracking, but $\tau(\R)$ touches itself.}
\label{fig:hole}
\end{figure}
Thus, after sufficiently many applications of the reduction operation, we obtain
$\sigma(\R)$, which has no back tracking, and no different edges where it touches itself.
As mentioned above, its left and right sides will touch paths of the form
$\sigma(\A)$, $\sigma(\B)$ for different {\A} and ${\B}$.

The fact that {\S} and {\s} are equivalent to reductions of $\texttt{Rvr}$ and $\texttt{rvR}$ is
because these start in the correct
parity squares (even and odd respectively), and right followed by reverse then right is the same as straight ahead, once the
backtracking is removed, as shown in Figure~\ref{fig:reducecases}.
This uses the assumption that all unit movements are by the same amount.
Since {\l} and {\r} are obtained from {\B} instead of {\A}, and the construction
of $\sigma(\B)$ is a reversing and complementing of
$\sigma(\A)$, it follows that
the turning directions for $\tau(\l)$ and $\tau(\r)$ are obtained by reversing those of
$\tau(\R)$ and $\tau(\L)$.

\subsection{Parity / case of letters of sequences}

The applications of the reduction algorithm results in a sequence of directions, in terms of
$\texttt{R,L,S}$, describing the boundary curve.
However, since application of $\tau$ depends on the parity of the starting square (which determines whether
the paths are around $\sigma(\A)$ or $\sigma(\B)$), we need to specify the
parity of the starting square of an instruction, which is indicated by case.
In Figure~\ref{fig:kindsofconnections}, it can be seen that the segment following a right or left turn starts from a
different parity, but following a straight ahead instruction does not change the
parity of the starting point of the next segment.
This is because a right turn will pass from, for example, square $(i,j)$ to $(i+1,j)$,
which has opposite parity,
whereas a straight ahead instruction will pass from square $(i,j)$ to one of $(i\pm 1,j\pm 1)$ all of which have the
same parity as the initial $(i,j)$ square.
Therefore, case of letters alternate, except when there is an {\S} or {\s}, in which case the
next letter has the same case.
Thus, the only thing to decide is the case of the first letter.
The parity of the squares is defined so that {\A} always starts in an even square.
So provided the first letter of $\sigma(\A)$ is {\A}, then the
first letters of $\tau(\R)$, $\tau(\L)$, $\tau(\S)$,
are capitals (denoting even start).
Similarly for the first letter of $\sigma(\B)$, if
$\sigma(\B)$ starts with \B, then paths starting from \B remain in odd squares, and in that
case the first letters of $\tau(\l)$, $\tau(\r)$, $\tau(\s)$ are lower case.
Since $\sigma(\A)$ determine $\sigma(\B)$, we can also write this in terms of $\sigma(\A)$, since the
last letter of $\sigma(\A)$ is the opposite of the first letter of $\sigma(\B)$ (calling $\A$, $\B$ opposite to each
other), and so
the first letters
of  $\tau(\l)$, $\tau(\r)$, $\tau(s)$ are lower case if and only if the last letter of $\sigma(\A)$
is \A.

\begin{figure}[H]
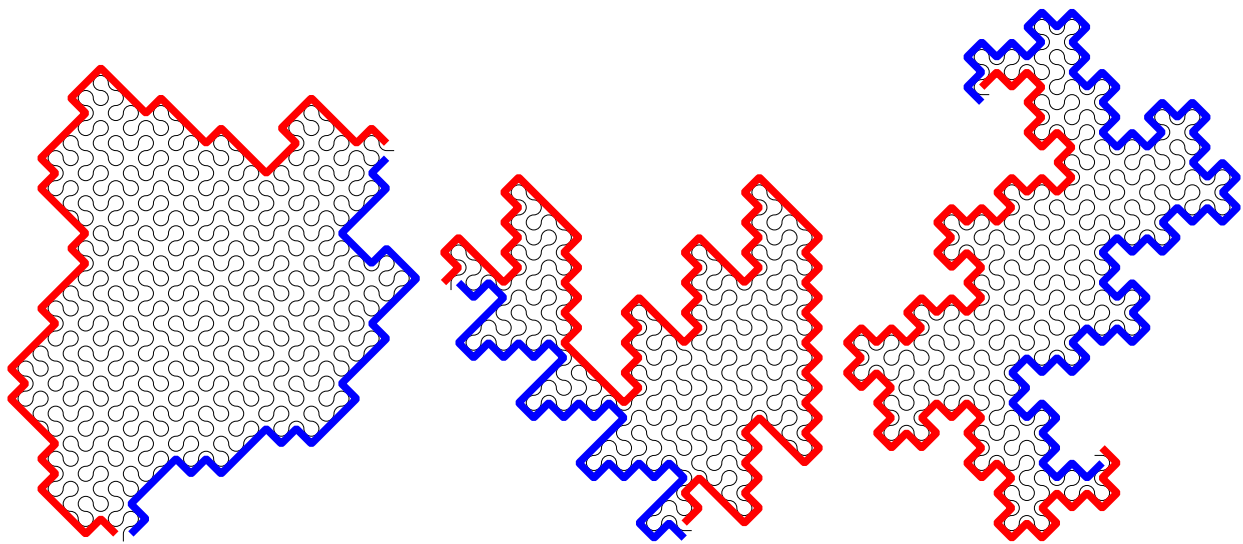

\input{exampleE10.tex}
\input{exampleE8.tex}
\input{exampleE5a.tex}
\caption{Curves and boundaries corresponding to the
  L-systems in Table~\ref{tab:lastexample}; iterations to
  level $3, 3, 4$ respectively.}
\label{fig:lastfig}
\end{figure}

The result of the above discussion is the following result, where
$\tau^\infty$  and $\sigma^\infty$ mean an appropriate limit, as for example in \cite{Edgar},
with the limit of $\tau^n$ chosen in a compatible way to $\sigma^n$.  For example, this can be achieved by
applying an appropriate isometric transformation $M_n$ to $\mathbb R^2$ such that the start and end points
of the paths $M_n(\sigma^n(\A))$, $M_n(\tau(\R))$, $M_n(\tau(\L))$ remain fixed.
\begin{theorem}
  The curves $\tau^\infty(\R)$ and
  $\tau^\infty(\L)$ form the boundary of the plane-filling curve
  $\sigma^\infty(\A)$.
\end{theorem}
Further examples can be found in Figure~\ref{fig:lastfig}, Table~\ref{tab:lastexample}, and in the implementation at \cite{verrill}.

\section*{Acknowledgments}

Many thanks to J.~Arndt for  encouragement and helpful comments on drafts.

{\small

\bibliographystyle{plain}
}

\begin{table}[H]
$$
  \begin{array}{|l|l|}
    \hline
1)&  \texttt{B+A+B-A-B-A+B+A+B-A}\\
2)&  \texttt{A+B-A-B+A+B+A-B}\\
3)& \texttt{A+B+A-B-A}\\
    \hline
  \end{array}
  \>\>
\begin{array}{|l||l|l|l|}
  \hline
  & (1) & (2) & (3)\\
      \hline
  \R &\texttt{rSRs}&\texttt{RlRrS}&\texttt{RrL}\\
  \r &\texttt{srL}&\texttt{SS}&\texttt{rRL}\\
  \L &\texttt{rLs}&\texttt{SS}&\texttt{RlL}\\
  \l &\texttt{slSL}&\texttt{SLlRl}&\texttt{rLl}\\
  \S &\texttt{rSSL}&\texttt{RlRl}&\texttt{Rsl}\\
  \s &\texttt{slRs}&\texttt{SLlRrS}&\texttt{rSL}\\
      \hline
\end{array}
$$
\caption{Example space filling L-systems and corresponding boundary
  L-systems}
\label{tab:lastexample}
\end{table}

\end{document}